\documentclass[a4paper,3p,11pt,sort&compress]{elsarticle}

\RequirePackage{snapshot}

\usepackage[Symbol]{upgreek}
\usepackage{mathptmx}
\usepackage[framed,numbered,autolinebreaks,useliterate]{mcode}
\usepackage[T1]{fontenc}
\usepackage[latin9]{inputenc}
\usepackage{amsmath}
\usepackage{amssymb}
\usepackage{graphicx}
\usepackage{setspace}
\usepackage{esint}
\PassOptionsToPackage{normalem}{ulem}
\usepackage{ulem}
\usepackage{marvosym}
\usepackage{multirow}
\usepackage{soul} 
\usepackage{physics}
\usepackage{xfrac}
\usepackage{dsfont}
\usepackage[hidelinks]{hyperref}
\usepackage{nomencl}
\usepackage{adjustbox}
\usepackage{mathtools}
\usepackage{nicefrac}
\usepackage[hang,flushmargin]{footmisc}
\usepackage{amsthm}
\newtheorem{theorem}{Theorem}

\DeclareMathAlphabet\mathbfcal{OMS}{cmsy}{b}{n}

\makeatletter
\newcommand{\shorteq}{%
  \settowidth{\@tempdima}{-}
  \resizebox{\@tempdima}{\height}{=}%
}
\makeatother

\DeclareMathSymbol{\shortminus}{\mathbin}{AMSa}{"39}
\newcommand{\shortm}{\negthinspace\shortminus\negthinspace}

\soulregister\cite7
\soulregister\ref7
\soulregister\pageref7

\def\*#1{\mathbf{#1}}

\usepackage{hyperref}

\makeatletter

\newcommand{\vt}[1]{\mathbf{#1}}                    
\newcommand{\vg}[1]{\bm{#1}}                        
\newcommand{\T}[0]{\mathrm{T}}                      
\newcommand{\I}[0]{\mathrm{i}\mkern1mu}                      
\newcommand{\e}[0]{\mathrm{e}}                      

\newcommand{\Ez}[0]{\mathbf{E}_0}
\newcommand{\Eo}[0]{\mathbf{E}_1}
\newcommand{\Et}[0]{\mathbf{E}_2}
\newcommand{\M}[0]{\mathbf{M}}
\newcommand{\Ezbar}[0]{\xoverline{\*E}_0}
\newcommand{\Eobar}[0]{\xoverline{\*E}_1}
\newcommand{\Etbar}[0]{\xoverline{\*E}_2}
\newcommand{\Mbar}[0]{\xoverline{\*M}}

\newcommand{\minusspace}[0]{\phantom{\shortm}} 			

\newcommand*{\eten}{\*C} 								

\newcommand*{\uv}{\xoverline{\*u}} 						
\newcommand*{\ua}{\*u} 									
\newcommand*{\uh}{\*u_h} 								
\newcommand*{\un}{\*u_n} 								
\newcommand*{\tfa}{\*v} 								
\newcommand*{\tfh}{\*v_h} 								
\newcommand*{\pf}{\xoverline{p}} 						
\newcommand*{\pa}{p} 						
\newcommand*{\ps}{p_s} 						
\newcommand*{\pn}{\*p_n} 								

\newcommand{\dotP}[0]{\scalebox{1}{\textbullet}}            

\newcommand{\egv}[0]{\vg{\upphi}}                  

\DeclareMathSymbol{\shortminus}{\mathbin}{AMSa}{"39}

\DeclareMathAlphabet\mathbfcal{OMS}{cmsy}{b}{n}      
\DeclareMathAlphabet\mathcal{OMS}{cmsy}{m}{n}      

\newcommand{\eq}{\begin{eqnarray}}
\newcommand{\nq}{\end{eqnarray}}
\newcommand{\dirvec}[1]{\mathbf{e}_{#1}}

\newcommand\ten[1]{\mathbf{#1}}

\usepackage{subfig} 
\usepackage{bm}
\usepackage{color}
\usepackage{epstopdf}

\newcommand*\xbar[1]{%
  \hbox{%
    \vbox{%
      \hrule height 0.5pt 
      \kern0.3ex
      \hbox{%
        \kern-0.1em
        \ensuremath{#1}%
        \kern-0.1em
      }%
    }%
  }%
}

\newsavebox\myboxA
\newsavebox\myboxB
\newlength\mylenA

\newcommand*\xoverline[2][0.75]{%
    \sbox{\myboxA}{$\m@th#2$}%
    \setbox\myboxB\null
    \ht\myboxB=\ht\myboxA%
    \dp\myboxB=\dp\myboxA%
    \wd\myboxB=#1\wd\myboxA
    \sbox\myboxB{$\m@th\overline{\copy\myboxB}$}
    \setlength\mylenA{\the\wd\myboxA}
    \addtolength\mylenA{-\the\wd\myboxB}%
    \ifdim\wd\myboxB<\wd\myboxA%
       \rlap{\hskip 0.5\mylenA\usebox\myboxB}{\usebox\myboxA}%
    \else
        \hskip -0.5\mylenA\rlap{\usebox\myboxA}{\hskip 0.5\mylenA\usebox\myboxB}%
    \fi}

\usepackage{tikz,pgfplots,grffile} 
\pgfplotsset{compat=newest}
\usetikzlibrary{calc}

\graphicspath{{figures/},{figures/sbdisp/}}

\definecolor{myred}{rgb}{0.6353,0.0784, 0.1843}
\definecolor{myblue}{rgb}{0, 0.2588, 0.6275}
\definecolor{mygreen}{rgb}{0, 0.5333, 0.2400}
\definecolor{mypurple}{rgb}{0.5882, 0.1176, 0.7059}
\definecolor{myorange}{rgb}{0.8431, 0.3333, 0.0980}
\definecolor{mygray}{rgb}{0.5020, 0.5020, 0.5020}

\journal{}

\newdefinition{rmk}{Remark}

\makeatother
\PassOptionsToPackage{numbers,sort&compress}{natbib}

\makenomenclature
\begin{document}
\makeatletter
\def\bm@pmb@#1{{%
      \setbox\tw@\hbox{$\m@th\mkern.25mu$}%
      \mathchoice
      \bm@pmb@@\displaystyle\@empty{#1}%
      \bm@pmb@@\textstyle\@empty{#1}%
      \bm@pmb@@\scriptstyle\defaultscriptratio{#1}%
      \bm@pmb@@\scriptscriptstyle\defaultscriptscriptratio{#1}}}
\makeatother

\title{Computation of leaky waves in layered structures coupled to unbounded media by exploiting multiparameter eigenvalue problems}

\author[cimne]{Hauke~Gravenkamp}
\ead{hgravenkamp@cimne.upc.edu}

\author[lj]{Bor~Plestenjak}
\ead{bor.plestenjak@fmf.uni-lj.si}

\author[espci]{Daniel~A.~Kiefer}
\ead{daniel.kiefer@espci.fr}

\author[kth]{Elias~Jarlebring\corref{cor1}}
\ead{eliasj@kth.se}

\address[cimne]{International Centre for Numerical Methods in Engineering (CIMNE), 08034 Barcelona, Spain}

\address[lj]{IMFM and Faculty of Mathematics and Physics, University of Ljubljana, Jadranska 19, SI-1000 Ljubljana, Slovenia}

\address[espci]{Institut Langevin, ESPCI Paris, Universit\'e PSL, CNRS, 75005 Paris, France}

\address[kth]{Department of Mathematics, NA group, KTH Royal Institute of Technology, 100 44 Stockholm, Sweden}

\cortext[cor1]{Corresponding author}

\begin{abstract}\noindent
  We present a semi-analytical approach to compute quasi-guided elastic wave modes in horizontally layered structures radiating into unbounded fluid or solid media. This problem is of relevance, e.g., for the simulation of guided ultrasound in embedded plate structures or seismic waves in soil layers over an elastic half-space. We employ a semi-analytical formulation to describe the layers, thus discretizing the thickness direction by means of finite elements. For a free layer, this technique leads to a well-known quadratic eigenvalue problem for the mode shapes and corresponding horizontal wavenumbers. Rigorously incorporating the coupling conditions to account for the adjacent half-spaces gives rise to additional terms that are nonlinear in the wavenumber. We show that the resulting nonlinear eigenvalue problem can be cast in the form of a multiparameter eigenvalue problem whose solutions represent the wave numbers in the plate and in the half-spaces. The multiparameter eigenvalue problem is solved numerically using recently developed algorithms.
\end{abstract}
\maketitle

\section{Introduction}\noindent
One of the classical problems in elastodynamics is the description of wave propagation in elastic layers that are in contact with another medium of infinite extent at one or both surfaces. We present a computational approach for this problem in the setting where the unbounded domains are assumed homogeneous and consist either of an acoustic fluid (i.e., governed by the scalar wave equation) or a linearly elastic isotropic solid. Such configurations are frequently encountered, particularly in two different fields of engineering: In soil dynamics and earthquake engineering, the analysis of wave propagation and vibration in layered soils and rock formations is of interest, where it can often be assumed that the layers radiate energy into a much larger domain of, e.g., water or soil \cite{Flitman1962,Kausel2022}. In this field of study, large structures and small frequencies are typically encountered. At very different spatial and temporal scales, understanding the behavior of high-frequency guided waves in thin-walled structures immersed or embedded in other media is essential in the context of ultrasonic nondestructive testing, material characterization, or sensor development \cite{Drinkwater2003,Pelat2011,Pistone2015}.

For the modeling of waves in layered media, numerous analytical, numerical, and semi-analytical approaches have been developed. Analytical and semi-analytical methods are appropriate to obtain dispersion curves of guided and quasi-guided waves, which will, hereinafter, be the focus of the present work. Wave propagation in individual layers with traction-free surfaces is well-understood, and implicit closed-form expressions for their dispersion relations are available \cite{Lamb1917c,Kausel2013b}. However, obtaining from these relations all solutions at a given frequency already requires sophisticated numerical root-finding algorithms. Extensions to layered systems and their coupling to unbounded domains have been achieved by means of the Transfer Matrix Method, Global Matrix Method, or Stiffness Matrix Method \cite{Knopoff1964,Nayfeh1991b,huberClassificationSolutionsGuided2023}. However, in such cases, the numerical difficulties in obtaining all solutions as roots of the dispersion relations are significant. For this reason, semi-analytical methods have been established as effective means of describing waves in plates and other structures of constant cross-section \cite{Dong1972}. These approaches involve a discretization of the structure's cross-section, which, in the case of layered media in two dimensions, reduces to a one-dimensional discretization along the layer's thickness. Typically, standard finite element spaces are employed to this end \cite{Kausel1977,Bartoli2006,Gravenkamp2012}, while variants based on, e.g., NURBS \cite{Gravenkamp2016} or spectral collocation \cite{Kiefer2022} exist. Such formulations naturally lead to a quadratic eigenvalue problem, whose solutions represent wavenumbers and discretized mode shapes and can be solved robustly and efficiently using conventional algorithms. To researchers in soil dynamics, this concept of semi-discretization is most famously known as the Thin Layer Method (TLM) \cite{Kausel2004}, while, in the context of ultrasonic testing, the term Semi-Analytical Finite Element Method (SAFE) \cite{Mazzotti2013b} is common. As a more general approach that has its origins in the TLM but is capable of describing more arbitrary bounded and unbounded star-convex domains, the Scaled Boundary Finite Element Method (SBFEM) has been developed and expanded since the 1990s \cite{Song1997,Krome2017}.

Notwithstanding the efficiency and robustness of semi-analytical methods for modeling waveguides of arbitrary cross-sections, incorporating an adjacent unbounded domain coupled to the structure of interest is not straightforward. Consequently, a variety of formulations have been suggested to address this challenge.
An obvious idea is to discretize a finite part of the unbounded domain large enough to allow waves to decay sufficiently such that reflections of the computational domain's boundary cannot interfere with the waveguide modes. These approaches require artificial damping in the unbounded domain, which is often achieved either by classical absorbing regions with artificial viscoelastic damping \cite{Castaings2008} or coordinate mapping techniques such as perfectly matched layers (PML) \cite{Treyssedea,Nguyen2013a,Berenger1994}, infinite elements \cite{Qi2010, Jia2011}, or similar approaches \cite{Georgiades2022}. Still, such discretizations are expensive and engender severe challenges in obtaining suitable damping properties and distinguishing correctly between modes in the waveguide and spurious modes in the damped medium. Hence, it is desirable to describe the effect of the unbounded domain on the waveguide solely by an interface condition, avoiding the additional discretization and the associated set of nonphysical parameters.

A particularly simple attempt at achieving this goal relies on approximating the effect of an unbounded domain by linear damping (\textit{dashpot} boundary condition) \cite{Gravenkamp2014b, Gravenkamp2015}. The implied assumption is that waves are radiated in the direction normal to the waveguide's surface into the surrounding medium, which can be a reasonable approximation in some cases, depending on the combination of material parameters and the frequency range. 
On the other hand, incorporating the exact boundary conditions representing the unbounded domain in a rigorous way gives rise to additional \textit{nonlinear} terms in the eigenvalue problem. For the plane geometries considered here, these nonlinear terms are typically of the form
\begin{equation}\label{eq:ky}
  \kappa_y = \sqrt{\kappa^2 - k^2},
\end{equation}
relating a vertical wavenumber~$\kappa_y$ in the unbounded medium to the common horizontal wavenumber~$k$ (i.e., the eigenvalue), with some constant~$\kappa$.
In the case of an adjacent fluid medium, an iterative solution scheme has been devised, which starts from the linear dashpot approximation and iteratively corrects each modal solution by updating the dependency of the fluid coupling on the wavenumber in the structure \cite{Gravenkamp2014c}. This approach is relatively efficient but not always robust, as different starting values may converge to the same solution, resulting in missing modes in the spectrum. Moreover, the trapped waves (quasi-Scholte modes) that appear as additional solutions compared to the free plate are not found. For the special case of a solid plate in contact with a fluid on one side or the same fluid on both sides, Kiefer et al.~\cite{Kiefer2019,kieferElastodynamicQuasiguidedWaves2022} showed that the eigenvalue problem can be linearized by a change of variables. This approach robustly finds the full spectrum but cannot be generalized to other cases. Hayashi et al.\ found solutions for the same special case by exploiting symmetries in the wave propagation that enabled the linearization of the problem \cite{Hayashi2014}. More recently, Tang et al.~\cite{tangStudyLeakyLamb2022} and Ducasse and Deschamps~\cite{Ducasse2022} independently presented linearization methods that consist in introducing a higher-dimensional state space. The latter studies consider two different fluids on both sides of the plate. 

For more general setups, including layers embedded in solid media or coupled to different materials on both sides, the solution becomes challenging due to the presence of different modes in the unbounded media. In this paper, we will rigorously derive the coupling conditions between elastic plates and unbounded fluid or solid media. To describe the plate, we employ a semi-analytical finite element method with high-order Lagrange elements (as discussed in \cite{Gravenkamp2012,Gravenkamp2014a}) in order to keep the size of the resulting eigenvalue problem to a minimum. An unbounded fluid domain is described by only one additional degree of freedom (DOF) representing the pressure at the interface. In the case of an unbounded solid, we use displacement amplitudes at the interface, hence adding two or three degrees of freedom (depending on whether we include shear-horizontal modes) for each unbounded domain.

Due to the incorporation of coupling conditions such as \eqref{eq:ky} and more advanced relations given in Sections~\ref{sec:fluid} and \ref{sec:solid}, the discretized problem is a nonlinear eigenvalue problem (NLEVP) in the sense of \cite{Mehrmann:2004:NLEVP}. The solution of such problems is an active research field, particularly in the scope of numerical linear algebra, with a plentitude of algorithms for general problems as well as approaches tailored to specialized structures \cite{Guttel2017}, including those that arise in the modeling of waveguides \cite{Jarlebring:2017:TIAR}. Algorithms and methods have also been included in HPC software packages such as SLEPc \cite{Campos:2019:NEP}.

For the particular type of nonlinearity encountered in the current work, we employ a rather different technique based on a relationship between NLEVPs and multiparameter eigenvalue problems, as formally defined below. The general connection was pointed out in \cite{Shao:2018:ALGNEP} and has been exploited to develop methods for such problems in \cite{Ringh_Jarlebring}. We will demonstrate that the nonlinear eigenvalue problem addressed here can be cast into the form of a multiparameter eigenvalue problem, in which the vertical wavenumbers of the form $\kappa_y$ are interpreted as additional unknown parameters related to the eigenvalue $k$ through (in the case of fluid loading) Eq.~\eqref{eq:ky}, or slightly different relations in the case of coupling to a solid medium.

Generally, in a multiparameter eigenvalue problem \cite{AtkinsonMingarelli}, we have $r$ equations of the form
\begin{equation}\label{eq:mupa}
  (\*A_{i0}+\lambda_1 \*A_{i1} + \cdots + \lambda_r \*A_{ir})\,\*x_i=\vt{0},\qquad i=1,\ldots,r,
\end{equation}
where $\*A_{ij}$ is an $n_i\times n_i$ matrix. Eigenvalues are tuples $(\lambda_1,\ldots,\lambda_r)$ for which nonzero $\*x_i$ exist such that all equations are satisfied. Problems of this type also appear, e.g., in separable boundary value problems \cite{PCH_SpecColl_JCP}. Recently, this formalism was applied to the computation of zero-group velocity (ZGV) points in waveguides \cite{Kiefer2022}. A standard numerical approach to computing all eigenvalues consists in constructing the associated system of generalized eigenvalue problems $\vg{\Delta}_i \*z = \lambda_i \vg{\Delta}_0 \*z$, $i=1,\ldots,r$, for details see \ref{sec:appx_mep}. The latter system has the same eigenvalues $(\lambda_1,\ldots,\lambda_r)$ as the original problem but can now be solved with conventional methods as shown in \cite{HKP_JD2EP}, since the $r$ equations decouple in the eigenvalues. If the problem is singular, the regular part is first extracted using a staircase-type method \cite{MP_Q2EP}. We believe that this mathematical solution procedure is a more general and abstract variant of the physics-based formulations described in \cite{tangStudyLeakyLamb2022,Ducasse2022}, but an in-depth analysis of the relation is outside the scope of this contribution.

In the following section, we briefly recap the semi-analytical formulation of a free plate and state the general form of the nonlinear eigenvalue problem resulting from including coupling conditions at the interfaces. Sections~\ref{sec:fluid} and \ref{sec:solid} explain the interface conditions for unbounded fluid and elastic media, respectively. The solution of the nonlinear eigenvalue problems is discussed in Section~\ref{sec:solution}, before we present some numerical examples (Section~\ref{sec:examples}) and a conclusion (Section~\ref{sec:conclusion}).

\section{Semi-analytical model of a layered medium}\label{sec:free}\noindent
The simple case of a free plate is revisited first, based on which we introduce the general effect of coupling to exterior half-spaces.

\subsection{Free plate}\label{sec:safe}\noindent

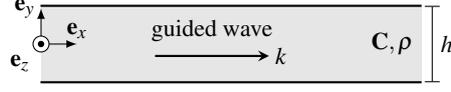
\begin{figure}[tb]\centering

\begin{tikzpicture}[>=stealth,
plate/.style={fill=black!10, minimum width=\l, minimum height=\h, inner sep=0pt},
]
\footnotesize
\pgfmathsetmacro\h{1 cm}
\pgfmathsetmacro\l{5 cm}
\pgfmathsetmacro\kxlen{1.5 cm}

\node[plate] (plate) at (0, 0) {};
\draw[thick] (plate.north west) -- (plate.north east); 
\draw[thick] (plate.south west) -- (plate.south east); 

\draw[Bar-Bar] ([xshift=4pt]plate.north east) -- node[right]{$h$} ([xshift=4pt]plate.south east);

\coordinate (O) at ([xshift=-1cm, yshift=-0.5em]plate.center);
\coordinate (Pkx) at ([xshift=\kxlen]O);
\draw[->, thick] (O) -- node[above, inner sep=2pt, yshift=3pt]{guided wave} (Pkx) node[right, inner sep=2pt]{$k$};

\node[left, yshift=0pt] at (plate.east) {$\eten, \rho$};

\begin{scope}[shift={(-0.5*5cm, 0)}]
    \draw[-latex, fill opacity = 1] (0, 0) -- (1.5em, 0) node[above, inner sep=2pt]{$\dirvec{x}$};
    \draw[-latex, fill opacity = 1] (0, 0) -- (0, 1.5em) node[left, inner sep=2pt]{$\dirvec{y}$};
    \draw[fill=white] (0,0) circle (3pt) node[below left, inner sep=4pt]{$\dirvec{z}$};
    \draw[fill=black] (0,0) circle (1pt);
\end{scope}


\end{tikzpicture}
  \caption{Plate with traction-free surfaces.}
  \label{fig:plate_geometry}
\end{figure}

\noindent
The basic semi-analytical formulation of waves in a free elastic plate is well-known \cite{Kausel1977,Bartoli2006,Gravenkamp2012} and only briefly recapitulated here. As depicted in Fig.~\ref{fig:plate_geometry}, let the plate of thickness $h$ occupy the domain $\Upomega$, defined by the open interval $\mathcal{I} = (\shortminus \tfrac{h}{2},\tfrac{h}{2})$ such that
\begin{equation*}
  \Upomega=\left\{\left. (x,\,y) \in \mathbb{R}^2\ \right|\ x \in  \mathbb{R},\ y \in \mathcal{I}\right\}.
\end{equation*}
Assuming the linearized momentum equation of elasticity, the displacement field $\uv(x,y,t)$ in the plate is governed by
\begin{equation} \label{eq:Navier-Cauchy}%
  \nabla \cdot (\eten : \nabla \uv) - \rho \partial_{tt}\uv  = \*0
\end{equation}
with the mass density $\rho$ and stiffness tensor $\eten$, and $\partial_{tt}$ denotes the second partial derivative with respect to time. Consider harmonic wave propagation along $x$ with a wavenumber $k$ and frequency $\omega$, i.e.,
\begin{equation}	\label{eq:plane_wave_ansatz}%
  \uv(x,y,t) = \ua(k, y, \omega)\, \e^{\I (k x - \omega t)}.
\end{equation}
Substituting Eq.~\eqref{eq:plane_wave_ansatz} into \eqref{eq:Navier-Cauchy}, dividing by $\e^{\I (k x - \omega t)}$, and introducing traction-free boundary conditions,\footnote{The typical traction-free von-Neumann boundary conditions are assumed here for notational convenience, while we may, alternatively, introduce Dirichlet conditions at the plate's surfaces.} the following one-dimensional boundary value problem (BVP) can be stated. For a given $\omega$, find $(k,\,\ua)$ such that
\begin{subequations}\label{eq:bvp}
  \begin{align}\label{eq:bvp1}
    \left[ (\I k)^2 \eten_{xx} + \I k\, \eten_{xy}  \partial_y + \I k\, \eten_{yx} \partial_y + \eten_{yy} \partial_y^2 + \omega^2 \rho \*I \right] \ua & = \*0 \qquad y \in \mathcal{I}     \\ \label{eq:bvp2}
    \left[ \I k \eten_{yx} + \eten_{yy}\partial_y \right] \ua                                                                                           & = \*0  \qquad y = \pm \tfrac{h}{2}
  \end{align}
\end{subequations}
with the 2$^\text{nd}$-order tensors $\eten_{ij} = \dirvec{i} \cdot \eten \cdot \dirvec{j}$ and the direction vectors $\dirvec{i}, i \in \{x, y, z\}$.\footnote{For an isotropic medium with the Lam\'e parameters $\lambda$ and $\mu$, we have
  \begin{equation*}\renewcommand{\arraystretch}{0.8} \arraycolsep=3pt
    \eten_{xx}=
    \left[\begin{array}{ccc}
        \lambda + 2\mu & 0   & 0   \\
        0              & \mu & 0   \\
        0              & 0   & \mu
      \end{array}\right], \qquad
    \eten_{yy}=
    \left[\begin{array}{ccc}
        \mu & 0              & 0   \\
        0   & \lambda + 2\mu & 0   \\
        0   & 0              & \mu
      \end{array}\right], \qquad
    \eten_{xy}= \eten_{yx}^\T =
    \left[\begin{array}{ccc}
        0       & \mu & 0 \\
        \lambda & 0   & 0 \\
        0       & 0   & 0
      \end{array}\right].
  \end{equation*}}
We denote by \mbox{$\mathbfcal{W} = {\*H}^1(\Omega)$} the function space of the continuous problem with the standard notation for Sobolev spaces. Multiplying Eq.~\eqref{eq:bvp1} by test functions $\tfa \in \mathbfcal{W}$ (i.e., taken from the same space as $\ua$), integrating over $\mathcal{I}$, and applying integration by parts to the third and fourth term yields the variational statement
\begin{equation}\label{eq:weakForm}
  -k^2 \big(\tfa, \eten_{xx} \ua\big) + \I k \big(\tfa, \eten_{xy} \partial_y \ua\big) - \I k \big(\partial_y\tfa, \*C_{yx} \ua\big) - \big(\partial_y\tfa, \eten_{yy} \partial_y\ua\big) + \omega^2 \big(\tfa,\rho \ua\big) = \*0,
\end{equation}
where $(\dotP,\dotP)$ denotes the $\mathrm{L}^2$-inner product in $\mathcal{I}$. The weak form \eqref{eq:weakForm} is discretized straightforwardly by defining adequate finite element spaces $\mathbfcal{W}_h\subset \mathbfcal{W}$. Replacing the trial and test functions in \eqref{eq:weakForm} with the discrete counterparts $\uh, \tfh \in \mathbfcal{W}_h$ defines the discretized BVP. Performing the integrations numerically, one obtains a square matrix for each of the terms in \eqref{eq:weakForm}, denoted as
\begin{equation}
  \left(-k^2\, \Ez + \I k\,\Eo-\Et+\omega^2\,\M\right)\un=\*0.\label{eq:EVPfreePlate}
\end{equation}
Here, $\un$ refers to the vector of coefficients (in our case, nodal displacements) of the finite element discretization.
Equation~\eqref{eq:EVPfreePlate} relates frequencies and horizontal wavenumbers and, hence, can be interpreted as a discretized dispersion relation. While this equation poses a quadratic two-parameter eigenvalue problem, it is commonly solved by choosing a set of frequencies and computing the corresponding wavenumbers after employing a suitable linearization.
\begin{rmk}
  In the case of undamped free plates or similarly simple scenarios, we can alternatively obtain dispersion curves by choosing a set of real wavenumbers and computing the corresponding frequencies. This can be computationally cheaper, as Eq.~\eqref{eq:EVPfreePlate} directly poses a linear eigenvalue problem for $\omega^2$. However, in the case of material damping or radiation into adjacent media, all wavenumbers of guided wave modes are complex; hence, we must solve Eq.~\eqref{eq:EVPfreePlate} at given frequencies, which we know to lie on a real interval.
\end{rmk} 
\begin{rmk}
  In finite-element-based formulations like the one outlined above, it is straightforward to incorporate any anisotropic material behavior of the plate into the stiffness tensor, see, e.g., \cite{Kausel1986a,Gravenkamp2012}. On the other hand, when dealing with elastic half-spaces (Section~\ref{sec:solid}), we will assume those to be isotropic.
\end{rmk}
\begin{rmk}
  We will see that the computational costs for solving the nonlinear eigenvalue problems discussed in the ensuing increase rapidly with the size of the finite-element matrices. Hence, it is crucial to choose an efficient discretization scheme. For this purpose, it is highly beneficial to use as few elements as possible (in our current application, one element per layer) of an order adequate to obtain reasonably accurate results in the desired frequency range.  As has been discussed in much detail previously, a suitable and rather common choice of trial and test functions is provided by Lagrange interpolation polynomials with nodes positioned at the Gauss-Lobatto points, allowing polynomial degrees of $100$ and more without numerical issues \cite{Gravenkamp2012,Gravenkamp2013a,Gravenkamp2014a}. Alternatives can be found in hierarchical shape functions or NURBS \cite{Gravenkamp2016}, leading to similar numerical properties. An overview (in the context of the general SBFEM) is provided in \cite{Gravenkamp2019}.
\end{rmk}

\subsection{General form of the nonlinear EVP}\noindent
In the following sections, we will see that incorporating the coupling to unbounded domains at the plate's top and bottom surfaces leads to a modified nonlinear eigenvalue problem of the general form
\begin{equation}
  \left(-k^2\, \Ezbar + \I k\,\Eobar-\Etbar+\omega^2\,\Mbar+\*R(k)\right)\egv=\vt{0}. \label{eq:NLEVPgeneral}
\end{equation}
The eigenvector $\egv$ contains the nodal displacements $\un$ of the plate and also up to six additional degrees of freedom representing the unbounded domains. The matrix function $\*R(k)$ involves non-polynomial terms in the eigenvalue $k$ and incorporates the interaction of the waveguide with a surrounding/adjacent medium depending on the physical setup.
Different approaches can be used to describe the waveguide as well as the surrounding medium. We formulate wave propagation inside elastic media in terms of displacements, while other formulations based on elastic potentials or mixed interpolations are possible. Waves in acoustic media can be defined in terms of the acoustic pressure or the velocity potential (we opt for the former).
We also note that it is generally possible to eliminate the degrees of freedom in the unbounded domain and only incorporate their effect on the waveguide through (nonlinear) Neumann boundary conditions as has been done, e.g., in \cite{Gravenkamp2014c}. However, that approach leads to a more complicated dependency of the additional terms in $\*R(k)$ on the eigenvalue $k$ and hinders the solution by the approaches that we suggest in the following sections.

\section{Elastic plate in contact with a fluid half-space}\label{sec:fluid}\noindent

\begin{figure}\centering
  \raisebox{1cm}{\subfloat{

\begin{tikzpicture}[>=stealth,
plate/.style={fill=black!10, minimum width=\l, minimum height=\h, inner sep=0pt},
]
\footnotesize
\pgfmathsetmacro\h{1 cm}
\pgfmathsetmacro\l{5 cm}
\pgfmathsetmacro\kxlen{1.5 cm}

\node[plate] (plate) at (0, 0) {};
\draw[thick,dashed] (plate.north west) -- (plate.north east); 
\draw[thick] (plate.south west) -- (plate.south east); 

\draw[Bar-Bar] ([xshift=4pt]plate.north east) -- node[right]{$h$} ([xshift=4pt]plate.south east);

\coordinate (O) at ([xshift=-1cm, yshift=-0.5em]plate.center);
\coordinate (Pkx) at ([xshift=\kxlen]O);
\draw[->, thick] (O) -- node[above, inner sep=2pt, yshift=3pt]{quasi-guided wave} (Pkx) node[right, inner sep=2pt]{$k$};

\node[left, yshift=0pt] at (plate.east) {$\eten, \rho$};
\node[above left, yshift=3pt] at (plate.north east) {$\tilde{c}, \tilde{\rho}$};

\begin{scope}[shift={(-0.5*5cm, 0)}]
    \draw[-latex, fill opacity = 1] (0, 0) -- (1.5em, 0) node[above, inner sep=2pt]{$\dirvec{x}$};
    \draw[-latex, fill opacity = 1] (0, 0) -- (0, 1.5em) node[left, inner sep=2pt]{$\dirvec{y}$};
    \draw[fill=white] (0,0) circle (3pt) node[below left, inner sep=4pt]{$\dirvec{z}$};
    \draw[fill=black] (0,0) circle (1pt);
\end{scope}

\coordinate (Oupp) at ($(O)+(0,0.5cm+0.5em+3pt)$);
\draw[->, thick]  (Oupp) -- node[pos=0.5,above left,outer sep=0pt,inner sep=2pt]{$\kappa$} +(1.5cm,0.5cm);
\end{tikzpicture}}}\qquad\qquad
  \subfloat{

\begin{tikzpicture}[>=stealth,
plate/.style={fill=black!10, minimum width=\l, minimum height=\h, inner sep=0pt},
]
\footnotesize
\pgfmathsetmacro\h{1 cm}
\pgfmathsetmacro\l{5 cm}
\pgfmathsetmacro\kxlen{1.5 cm}

\node[plate] (plate) at (0, 0) {};
\draw[thick,dashed] (plate.north west) -- (plate.north east); 
\draw[thick,dashed] (plate.south west) -- (plate.south east); 

\draw[Bar-Bar] ([xshift=4pt]plate.north east) -- node[right]{$h$} ([xshift=4pt]plate.south east);

\coordinate (O) at ([xshift=-1cm, yshift=-0.5em]plate.center);
\coordinate (Pkx) at ([xshift=\kxlen]O);
\draw[->, thick] (O) -- node[above, inner sep=2pt, yshift=3pt]{quasi-guided wave} (Pkx) node[right, inner sep=2pt]{$k$};

\node[left, yshift=0pt] at (plate.east) {$\eten, \rho$};
\node[above left, yshift=3pt] at (plate.north east) {$\tilde{c}_1, \tilde{\rho}_1$};
\node[below left, yshift=-3pt] at (plate.south east) {$\tilde{c}_2, \tilde{\rho}_2$};

\begin{scope}[shift={(-0.5*5cm, 0)}]
    \draw[-latex, fill opacity = 1] (0, 0) -- (1.5em, 0) node[above, inner sep=2pt]{$\dirvec{x}$};
    \draw[-latex, fill opacity = 1] (0, 0) -- (0, 1.5em) node[left, inner sep=2pt]{$\dirvec{y}$};
    \draw[fill=white] (0,0) circle (3pt) node[below left, inner sep=4pt]{$\dirvec{z}$};
    \draw[fill=black] (0,0) circle (1pt);
\end{scope}

\coordinate (Oupp) at ($(O)+(0,0.5cm+0.5em+3pt)$);
\draw[->, thick]  (Oupp) -- node[pos=0.5,above left,outer sep=0pt,inner sep=2pt]{$\kappa_1$} +(1.5cm,0.5cm);

\coordinate (Obot) at ($(O)-(0,0.5cm-0.5em+3pt)$);
\draw[->, thick]  (Obot) -- node[pos=0.5,below left,outer sep=0pt,inner sep=2pt]{$\kappa_2$} +(1.5cm,-0.9cm);

\end{tikzpicture}}
  \caption{Plate in contact with fluid half-spaces at one or both of its surfaces. \label{fig:plate_geometry_fluid}}
\end{figure}
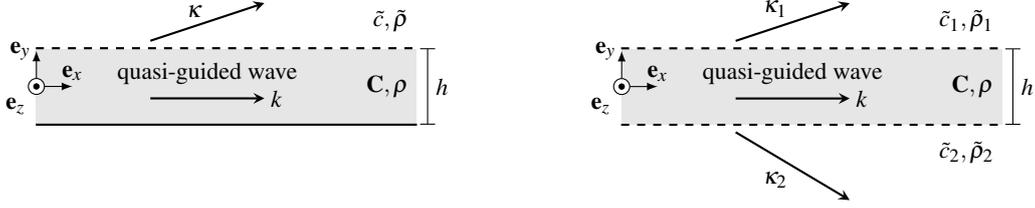

\noindent
Let us begin by assuming that the elastic plate, as described in the previous section, is now coupled to an infinite acoustic (inviscid) fluid at one of its surfaces located at $y=y_s$, see Fig.~\ref{fig:plate_geometry_fluid}. Hence, the acoustic pressure $\pf$ in the fluid satisfies the Helmholtz equation
\begin{equation}
  \Delta \pf + \kappa^2 \pf=0,
\end{equation}
where $\kappa$ denotes the wavenumber in the fluid. As we are dealing with a plane surface radiating into an infinite medium, we postulate plane pressure waves~\cite{viggenModellingAcousticRadiation2023} of amplitude $\pa$ propagating in the fluid domain, i.e.,
\begin{equation}\label{eq:acousticPlaneWave}
  \pf(x,y,t) = \pa\, \e^{\I(\kappa_x\,x + \kappa_y\,y-\omega\, t)},
\end{equation}
with $\kappa_x$ and $\kappa_y$ being the horizontal and vertical components of the complex-valued wave vector in the fluid. We note that the horizontal component matches the wavenumber of the guided waves inside the plate ($\kappa_x = k$), and the vertical wavenumber is given as
\begin{equation}\label{eq:verticalK}
  \kappa_y = \pm\sqrt{\kappa^2 - k^2}.
\end{equation}
Defining $\ps = \pa\, \e^{\I \kappa_y\, y_s}$, the fluid pressure along the interface reads
\begin{equation}\label{eq:acousticXT}
  \pf(x,y_s,t) = \ps\, \e^{\I(k x - \omega t)}.
\end{equation}
The coupling conditions at the interface are as follows. The acoustic pressure induces a traction $\vg{\tau}$ on the plate surface, i.e.,
\begin{equation}
  \vg{\tau}(y\,\shorteq\,y_s) = \*n\cdot\vg{\sigma}(y\,\shorteq\,y_s) = \ps\,\*n
\end{equation}
with the unit outward normal vector $\*n$ and the stress tensor $\vg{\sigma}$. Since, in our case, we have either $\*n = \*{e}_y$ or $\*n = -\*{e}_y$, this condition simplifies to
\begin{equation}\label{eq:couplingAcousticTraction}
  \tau_y(y\,\shorteq\,y_s) = \pm\ps,
\end{equation}
where the positive/negative sign corresponds to coupling at the upper/lower surface, respectively.
The second coupling condition states that the normal derivative of the acoustic pressure is related to the acceleration at the plate's surface by
\begin{equation}\label{eq:pNormalDerivative}
  \left.\*n\cdot\nabla \pf\right|_{y\,\shorteq\,y_s} = -\tilde{\rho}\,\*n\cdot \partial_t^2 \uv = \omega^2\tilde{\rho}\,\*n\cdot \uv
\end{equation}
with the fluid's density $\tilde{\rho}$. Hence, considering Eqs.~\eqref{eq:acousticPlaneWave} and \eqref{eq:verticalK}, we obtain
\begin{equation}\label{eq:couplingAcousticPressure}
  \pm \I \sqrt{\kappa^2 - k^2}\ps
  = \omega^2\tilde{\rho}\,u_y(y\,\shorteq\,y_s).
\end{equation}
Two possibilities exist for including the coupling conditions in the waveguide model. The first one, described in detail in \cite{Gravenkamp2014c}, consists in substituting  Eq.~\eqref{eq:couplingAcousticPressure} into \eqref{eq:couplingAcousticTraction} to obtain an expression for the traction while eliminating the acoustic pressure:
\begin{equation}\label{eq:couplingAcousticTractionNP}
  \tau_y(y\,\shorteq\,y_s)  = \mp\frac{\I  \omega^2\tilde{\rho}}{\sqrt{\kappa^2 - k^2}} \,u_y(y\,\shorteq\,y_s).
\end{equation}
The above equation can be considered as a Neumann boundary condition and directly integrated into the waveguide eigenproblem. However, the resulting nonlinear terms of the form $1/\sqrt{\kappa^2 - k^2}$, which are singular at $\kappa = k$, are somewhat difficult to address. Thus, we will follow the formulation similar to \cite{Kiefer2019}\footnote{In \cite{Kiefer2019}, the fluid was described in terms of the velocity potential rather than the pressure, and the discretization of the plate was achieved using spectral collocation rather than finite elements. However, these differences are of little relevance for the current discussion and the solution procedure proposed in this paper.} and introduce additional degrees of freedom $\pn$ representing the pressure amplitude(s) at the lower and upper surface. The eigenvectors are extended accordingly as
\begin{equation}
  \egv = \begin{bmatrix}
    \un \\ \pn
  \end{bmatrix}.
\end{equation}
In this approach, the coupling conditions are trivially incorporated into the nonlinear eigenvalue problem \eqref{eq:NLEVPgeneral} (extending the involved matrices to account for the additional DOFs), giving rise to a coupling matrix of the form
\begin{equation}\label{ex:nonpoly1}
  \*R(k) = \I\kappa_{y,1}\, \*R_1 + \I\kappa_{y,2}\, \*R_2,
\end{equation}
with 
\begin{equation}
  \kappa_{y,1} = \sqrt{\kappa_1^2 - k^2}, \quad  \kappa_{y,2} = \sqrt{\kappa_2^2 - k^2}.
\end{equation}
Here, the indices $1,2$ refer to the fluids coupled to the lower and upper surface. The additional matrix entries due to the coupling can be read off from Eqs.\ \eqref{eq:couplingAcousticTraction} and \eqref{eq:couplingAcousticPressure}.
Specifically, if we denote by $s_j$ the index of the DOF corresponding to the vertical displacement at interface $j=1,2$ and by $f_j$ the index of the pressure-DOF, these entries are
\begin{align*}
  \Etbar^{s_j,f_j} & = \pm 1            \\
  \Mbar^{f_j,s_j}  & = \pm \tilde{\rho} \\
  \*R_j^{f_j,f_j}  & = 1.
\end{align*}
Overall, this means that we consider \eqref{eq:couplingAcousticPressure} to be an additional equation incorporated into $\Mbar$ and $\*R_j$, while \eqref{eq:couplingAcousticTraction} is included in $\Etbar$ and represents the Neumann boundary condition imposed on the plate.
\begin{rmk}
  We may stress again that, in the case of $\kappa_1 = \kappa_2$ (i.e., a plate in contact with the same fluid medium at the lower and upper surface), the resulting eigenvalue problem can be linearized as discussed in \cite{Kiefer2019}. While we will study such an example in our numerical experiments to verify the implementation, the main objective of the approach presented here is to address more complex problems, particularly those involving coupling to solid media or different materials (elastic or acoustic) at both surfaces.
\end{rmk}

\section{Elastic plate embedded in an elastic medium}\label{sec:solid}\noindent%
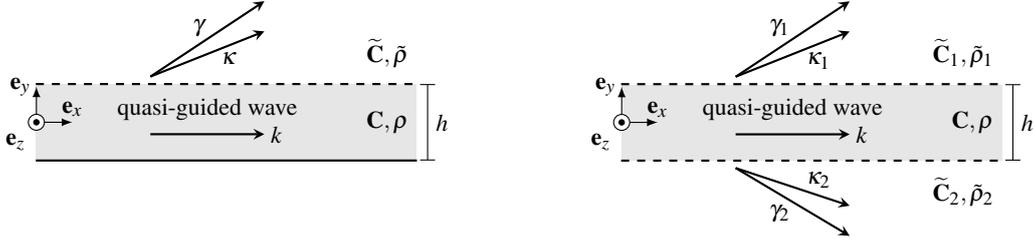
\begin{figure}\centering
  \raisebox{1cm}{\subfloat{

\begin{tikzpicture}[>=stealth,
plate/.style={fill=black!10, minimum width=\l, minimum height=\h, inner sep=0pt},
]
\footnotesize
\pgfmathsetmacro\h{1 cm}
\pgfmathsetmacro\l{5 cm}
\pgfmathsetmacro\kxlen{1.5 cm}

\node[plate] (plate) at (0, 0) {};
\draw[thick, dashed] (plate.north west) -- (plate.north east); 
\draw[thick] (plate.south west) -- (plate.south east); 

\draw[Bar-Bar] ([xshift=4pt]plate.north east) -- node[right]{$h$} ([xshift=4pt]plate.south east);

\coordinate (O) at ([xshift=-1cm, yshift=-0.5em]plate.center);
\coordinate (Pkx) at ([xshift=\kxlen]O);
\draw[->, thick] (O) -- node[above, inner sep=2pt, yshift=3pt]{quasi-guided wave} (Pkx) node[right, inner sep=2pt]{$k$};

\node[left, yshift=0pt] at (plate.east) {$\ten{C}, \rho$};
\node[above left, yshift=3pt] at (plate.north east) {$\widetilde{\ten{C}}, \tilde{\rho}$};

\begin{scope}[shift={(-0.5*5cm, 0)}]
    \draw[-latex, fill opacity = 1] (0, 0) -- (1.5em, 0) node[above, inner sep=2pt]{$\dirvec{x}$};
    \draw[-latex, fill opacity = 1] (0, 0) -- (0, 1.5em) node[left, inner sep=2pt]{$\dirvec{y}$};
    \draw[fill=white] (0,0) circle (3pt) node[below left, inner sep=4pt]{$\dirvec{z}$};
    \draw[fill=black] (0,0) circle (1pt);
\end{scope}

\coordinate (Oupp) at ($(O)+(0,0.5cm+0.5em+3pt)$);
\draw[->, thick]  (Oupp) -- node[pos=0.6,below right,outer sep=0pt,inner sep=1pt]{$\kappa$} +(1.5cm,0.6cm);
\draw[->, thick]  (Oupp) -- node[above left,outer sep=0pt,inner sep=1pt]{$\gamma$} +(1.5cm,1cm);


\end{tikzpicture}}}\qquad\qquad
  \subfloat{

\begin{tikzpicture}[>=stealth,
plate/.style={fill=black!10, minimum width=\l, minimum height=\h, inner sep=0pt},
]
\footnotesize
\pgfmathsetmacro\h{1 cm}
\pgfmathsetmacro\l{5 cm}
\pgfmathsetmacro\kxlen{1.5 cm}

\node[plate] (plate) at (0, 0) {};
\draw[thick, dashed] (plate.north west) -- (plate.north east); 
\draw[thick,dashed] (plate.south west) -- (plate.south east); 

\draw[Bar-Bar] ([xshift=4pt]plate.north east) -- node[right]{$h$} ([xshift=4pt]plate.south east);

\coordinate (O) at ([xshift=-1cm, yshift=-0.5em]plate.center);
\coordinate (Pkx) at ([xshift=\kxlen]O);
\draw[->, thick] (O) -- node[above, inner sep=2pt, yshift=3pt]{quasi-guided wave} (Pkx) node[right, inner sep=2pt]{$k$};

\node[left, yshift=0pt] at (plate.east) {$\ten{C}, \rho$};
\node[above left, yshift=3pt] at (plate.north east) {$\widetilde{\ten{C}}_1, \tilde{\rho}_1$};
\node[below left, yshift=-3pt] at (plate.south east) {$\widetilde{\ten{C}}_2, \tilde{\rho}_2$};

\begin{scope}[shift={(-0.5*5cm, 0)}]
    \draw[-latex, fill opacity = 1] (0, 0) -- (1.5em, 0) node[above, inner sep=2pt]{$\dirvec{x}$};
    \draw[-latex, fill opacity = 1] (0, 0) -- (0, 1.5em) node[left, inner sep=2pt]{$\dirvec{y}$};
    \draw[fill=white] (0,0) circle (3pt) node[below left, inner sep=4pt]{$\dirvec{z}$};
    \draw[fill=black] (0,0) circle (1pt);
\end{scope}

\coordinate (Oupp) at ($(O)+(0,0.5cm+0.5em+3pt)$);
\draw[->, thick]  (Oupp) -- node[pos=0.6,below right,outer sep=0pt,inner sep=1pt]{$\kappa_1$} +(1.5cm,0.6cm);
\draw[->, thick]  (Oupp) -- node[above left,outer sep=0pt,inner sep=1pt]{$\gamma_1$} +(1.5cm,1cm);

\coordinate (Obot) at ($(O)-(0,0.5cm-0.5em+3pt)$);
\draw[->, thick]  (Obot) -- node[pos=0.6,above right,outer sep=0pt,inner sep=1pt]{$\kappa_2$} +(1.5cm,-0.5cm);
\draw[->, thick]  (Obot) -- node[below left,outer sep=0pt,inner sep=1pt]{$\gamma_2$} +(1.5cm,-0.9cm);

\end{tikzpicture}}
  \caption{Plate in contact with solid half-spaces at one or both of its surfaces. \label{fig:plate_geometry_embedded}}
\end{figure}%
We now consider the situation in which the elastic plate is in contact with an elastic half-space or embedded in an elastic medium. For the beginning, let us again assume the presence of only one half-space coupled to the plate at $y=y_s$; see Fig.~\ref{fig:plate_geometry_embedded}. The extension to two half-spaces of potentially different materials coupled to the lower and upper surfaces is straightforward.
The fundamental solution of waves propagating in an isotropic horizontally stratified medium (assuming plane strain) can be found in the literature \cite{Rose1999,Kausel2006b} and is only briefly reproduced here to fix notation. Again, we consider the waveguide to be aligned with the $x$-axis; thus, the wavenumber $k$ of the guided wave modes matches the $x$-component of the wave vector in the surrounding medium.
The displacements in the unbounded solid are a superposition of three plane waves, namely one longitudinal and two transverse waves. The latter are denoted as shear-vertical and shear-horizontal waves, respectively.
The wavenumbers of longitudinal and transverse waves are denoted as $\kappa$ and $\gamma$, and the respective vertical components of the wave vector are
\begin{subequations}\label{eq:alphaBeta}
  \begin{align}
    \kappa_y & =\pm\sqrt{\kappa^2-k^2},  \\
    \gamma_y & = \pm\sqrt{\gamma^2-k^2},
  \end{align}
\end{subequations}
leading to nonlinearities similar to the case of fluid loading discussed before. 
Accounting for these observations, we write the displacements~ $\uv = (\xoverline{u}_x,\, \xoverline{u}_y,\, \xoverline{u}_z )^\T$  in the unbounded solid as
\begin{equation}\label{eq:displacementPotential}
  \uv(x,y,t) =
  \nabla \varphi +
  \nabla \times \bm{\uppsi} +
  \partial_x  \bm{\chi}
  , \qquad y \notin \mathcal{I},
\end{equation}
where $\varphi$ is the scalar potential of the longitudinal wave, $\bm{\uppsi}$ represents the vector potential of the shear-vertical wave, and $\bm{\chi}$ describes the shear-horizontal wave. The potentials are assumed as
\begin{subequations}\label{eq:potentialsAnsatz}
  \begin{align}
    \varphi     & = a\, \e^{\I(k x + \kappa_y y -\omega t)},             \\
    \bm{\uppsi} & = b\, \dirvec{z} \e^{\I(k x + \gamma_y  y -\omega t)}, \\
    \bm{\chi}   & = c\, \dirvec{z}\, \e^{\I(k x + \gamma_y y -\omega t)}
  \end{align}
\end{subequations}
with unknown wave amplitudes $a,b,c$.
Substituting Eqs.~\eqref{eq:potentialsAnsatz} into \eqref{eq:displacementPotential} yields
\begin{equation}\label{eq:displacements_exterior}
  \uv(x,y,t) = \left[
    (\I k\, \dirvec{x} + \I \kappa_y\, \dirvec{y})\, a\, \e^{\I \kappa_y y}  +
    (\I \gamma_y\, \dirvec{x} - \I k\, \dirvec{y})\, b\, \e^{\I \gamma_y y} +
    \I k\, \dirvec{z}\, c\, \e^{\I \gamma_y y} \right] \e^{\I(k x -\omega t)},
\end{equation}
which we may write compactly as
\begin{equation}\label{eq:u_general}
  \uv(x,y,t)= \mathbfcal{A}(k)\, \e^{\I\mathbfcal{D} y}\, \e^{\I (kx-\omega t)}\,\vt{c}
\end{equation}
with
\begin{equation*}
  \mathbfcal{A}(k)=\left[\begin{array}{ccc}
      \I k        & \phantom{\shortm\,}\I\gamma_y & 0    \\
      \I \kappa_y & \shortm\,\I k                 & 0    \\
      0           & 0                             & \I k
    \end{array}\right],
\end{equation*}
\begin{equation*}\label{eq:Ey}
  \mathbfcal{D}=\operatorname{diag}[\kappa_y, \gamma_y, \gamma_y],\qquad
  \*c = (a,b,c)^\T.
\end{equation*}
Compared to the formulation used for describing a surrounding fluid, we have now introduced the unknowns $\*c$ rather than the acoustic pressure.
We will also require the amplitudes of tractions acting on a plane parallel to the $x$-$z$-plane, cf.\ Section~\ref{sec:safe}:
\begin{equation}\label{eq:tractions}
  \xbar{\vg{\tau}}\,(x,y) = \widetilde{\eten}_{yx}\, \partial_x \uv(x,y) + \widetilde{\eten}_{yy}\, \partial_y \uv(x,y).
\end{equation}
Substituting Eqs.~\eqref{eq:u_general} into \eqref{eq:tractions} yields
\begin{equation}\label{eq:tau_halfspace}
  \xbar{\vg{\tau}}\,(x,y)= \I k\, \widetilde{\eten}_{yx}\, \mathbfcal{A}\, \e^{\I\mathbfcal{D} y}\, \e^{\I (kx-\omega t)}\vt{c}  +  \I\widetilde{\eten}_{yy}\, \mathbfcal{A}\, \mathbfcal{D}\, \e^{\I\mathbfcal{D} y}\, \e^{\I (kx-\omega t)}\vt{c}
\end{equation}
To define the coupling conditions, we evaluate the displacement amplitudes at the interface $y=y_s$ based on Eq.~\eqref{eq:u_general} and obtain
\begin{equation}
  \ua({y\,\shorteq\,y_s}) = \mathbfcal{A}\,\vt{c}_s \label{eq:cAu}
\end{equation}
with
\begin{equation}
  \vt{c}_s = \e^{\I\mathbfcal{D} y_s}\,\*c.
\end{equation}
Equation~\eqref{eq:cAu} is expanded as
\begin{equation}\label{eq:solidCouplingU}
  \ua({y\,\shorteq\,y_s})  =
  (\I k\mathbfcal{A}_0 +\I \kappa_y  \mathbfcal{A}_1 +\I \gamma_y \mathbfcal{A}_2)\*c_s
\end{equation}
with
\begin{equation}
  \mathbfcal{A}_0 =
  \left[\begin{array}{rrr}
      1 & 0         & 0 \\
      0 & \shortm 1 & 0 \\
      0 & 0         & 1
    \end{array}\right], \qquad
  {\mathbfcal{A}}_1 =
  \left[\begin{array}{rrr}
      0 & 0 & 0 \\
      1 & 0 & 0 \\
      0 & 0 & 0
    \end{array}\right], \qquad
  \mathbfcal{A}_2 =
  \left[\begin{array}{rrr}
      0 & 1 & 0 \\
      0 & 0 & 0 \\
      0 & 0 & 0
    \end{array}\right].
\end{equation}
The second coupling condition is obtained by evaluating the traction amplitudes at the interface:
\begin{align*}
  \vg{\tau}({y\,\shorteq\,y_s})
   & = \I k \, \widetilde{\eten}_{xy}\, \mathbfcal{A}\,\*c_s  +  \I\widetilde{\eten}_{yy}\, \mathbfcal{A}\, \mathbfcal{D}\,\*c_s \\
   & = -k\widetilde{\eten}_{xy}\, (k\mathbfcal{A}_0 + \kappa_y  \mathbfcal{A}_1 + \gamma_y \mathbfcal{A}_2)\,\*c_s
  -  \widetilde{\eten}_{yy}\, ( k\mathbfcal{A}_0 +  \kappa_y  \mathbfcal{A}_1 + \gamma_y \mathbfcal{A}_2)\, \mathbfcal{D}\,\*c_s,
\end{align*}
where we have substituted Eqs.~\eqref{eq:Dmatrices} and \eqref{eq:solidCouplingU}.
We further introduce
\begin{equation}\label{eq:Dmatrices}
  \mathbfcal{D}_1 = \operatorname{diag}[1,0,0],\quad
  \mathbfcal{D}_2 = \operatorname{diag}[0,1,1]
\end{equation}
such that $\mathbfcal{D} = \kappa_y \mathbfcal{D}_1 + \gamma_y \mathbfcal{D}_2$ and note that
\begin{equation*}
  \mathbfcal{A}_1\,\mathbfcal{D}_2=\mathbfcal{A}_2\,\mathbfcal{D}_1 = \*0,\quad
  \mathbfcal{A}_1\,\mathbfcal{D}_1 = \mathbfcal{A}_1, \quad
  \mathbfcal{A}_2\,\mathbfcal{D}_2 = \mathbfcal{A}_2, \quad
  \mathbfcal{A}_0\,\mathbfcal{D}_1 = \mathbfcal{D}_1, \quad
  \mathbfcal{A}_0\,\mathbfcal{D}_2 = -\mathbfcal{D}_2.
\end{equation*}
Hence, the above expression for the traction is simplified as
\begin{equation}\label{eq:tractionSolid}
  \vg{\tau}({y\,\shorteq\,y_s})  =(k^2 \mathbfcal{R}_0  + k\kappa_y\, \mathbfcal{R}_1 + k\gamma_y\, \mathbfcal{R}_2 + \kappa^2 \mathbfcal{R}_3 + \gamma^2 \mathbfcal{R}_4)\,\*c_s
\end{equation}
with
\begin{align*}
  \mathbfcal{R}_0 & = -(
  \widetilde{\eten}_{xy}\mathbfcal{A}_0
  -\widetilde{\eten}_{yy}\mathbfcal{A}_1
  - \widetilde{\eten}_{yy}\mathbfcal{A}_2),                   \\
  \mathbfcal{R}_1 & = -(
  \widetilde{\eten}_{xy} \mathbfcal{A}_1
  +\widetilde{\eten}_{yy} \mathbfcal{D}_1),                   \\
  \mathbfcal{R}_2 & = - (
  \widetilde{\eten}_{xy} \mathbfcal{A}_2
  - \widetilde{\eten}_{yy} \mathbfcal{D}_2),                  \\
  \mathbfcal{R}_3 & = -\widetilde{\eten}_{yy}\mathbfcal{A}_1, \\
  \mathbfcal{R}_4 & = -\widetilde{\eten}_{yy}\mathbfcal{A}_2.
\end{align*}
To describe the plate coupled to an infinite solid medium, we introduce the conditions defined by Eqs.~\eqref{eq:solidCouplingU} and \eqref{eq:tractionSolid} (for two different media at the two surfaces if necessary) into the nonlinear eigenvalue problem. To this end, we multiply Eq.~\eqref{eq:solidCouplingU} by $\I k$;
hence, all coupling terms are either proportional to $k\kappa_y$ or $k\gamma_y$ or can be integrated into the existing matrices $\Ezbar$, $\Eobar$, $\Etbar$. 
The eigenvector contains the displacements in the plate and the modal amplitudes in the unbounded media, which we abbreviate as $\*c_n$:
\begin{equation}
  \egv = \begin{bmatrix}
    \un \\ \*c_n
  \end{bmatrix}.
\end{equation}
In the general case of two half-spaces consisting of different materials, the additional matrix function in the nonlinear eigenvalue problem is of the form
\begin{equation}\label{eq:nonpoly2}
  \*R(k) = k\kappa_{y,1}\, \*R_{1,1} + k\gamma_{y,1}\, \*R_{2,1} + k\kappa_{y,2}\, \*R_{1,2} + k\gamma_{y,2}\, \*R_{2,2}.
\end{equation}
Denoting by $p_j$ and $q_j$ the DOFs corresponding to the plate and unbounded domain at interface~$j$, the components introduced by the coupling conditions are
\begin{align*}
   & \Ezbar^{q_j,q_j} =  \mathbfcal{A}_0, \qquad
   &                                                      & \Ezbar^{p_j,q_j} =  \mp \mathbfcal{R}_{0,j}, \qquad
   &                                                      & \Eobar^{q_j,p_j} =  -\*I, \qquad
   &                                                      & \Etbar^{p_j,q_j} =  \mp\kappa_{j}^2\,\mathbfcal{R}_{3,j} \mp\gamma_{j}^2\,\mathbfcal{R}_{4,j}, \\
   & \*R_{1,j}^{p_j,q_j} = \pm\mathbfcal{R}_{1,j}, \qquad
   &                                                      & \*R_{2,j}^{p_j,q_j} = \pm\mathbfcal{R}_{2,j}, \qquad
   &                                                      & \*R_{1,j}^{q_j,q_j} = -\mathbfcal{A}_1, \qquad
   &                                                      & \*R_{2,j}^{q_j,q_j} = -\mathbfcal{A}_2.
\end{align*}

\section{Solution procedure}\label{sec:solution}\noindent
This section summarizes the approach we employ for solving the nonlinear eigenvalue problems derived before. The key is to realize that, in the particular cases encountered here, we can rewrite the NLEVP as a linear multiparameter eigenvalue problem. In particular, this is achieved by generalizing the transformation in \cite[Lemma~2.5]{Ringh_Jarlebring}, which relies on introducing auxiliary variables. For the so-obtained multiparameter eigenvalue problem, established algorithms exist, and their implementation is publicly available \cite{multipareig_2023}. Hence, we will not delve into the numerical details of the computation. For the sake of tangibility, we present, in this section, a rather general and clear form of the multiparameter eigenvalue problem, albeit it is not always the most efficient version. Some comments and suggestions for improving efficiency are given in the appendix.
The following theorem illustrates how the problem can be cast into a multiparameter 
eigenvalue problem.
\begin{theorem} If $(k,\egv)$ is a solution to \eqref{eq:NLEVPgeneral} with $\*R(k)$ given by \eqref{ex:nonpoly1}, then the parameters $\I k,\,\I\kappa_{y,1},\,\I\kappa_{y,2},\,\xi_0$ form a solution to the four-parameter eigenvalue problem
\begin{subequations}\label{eq:4par}
\begin{align}
  \left(-\Etbar+\omega^2\,\Mbar + \I k\,\Eobar + \I \kappa_{y,1}\*R_1 + \I \kappa_{y,2} \*R_2 + \xi_0\, \Ezbar \right)\egv & =\vt{0}\label{eq:4par_main_a}\\
  \left(\begin{bmatrix*}[r] 0 & \shortm\kappa_1^2 \cr 1 & 0 \end{bmatrix*}
  + \I\kappa_{y,1} \begin{bmatrix*}[r] 1 & 0 \cr 0 & 1 \end{bmatrix*}
  + \xi_0 \begin{bmatrix*}[r] 0 & \shortm 1 \cr 0 & 0 \end{bmatrix*}\right) \*x_1                                                       & = \vt{0} \label{eq:4par_main_b}\\
  \left(\begin{bmatrix*}[r] 0 & \shortm\kappa_2^2 \cr 1 & 0 \end{bmatrix*}
  + \I\kappa_{y,2} \begin{bmatrix*}[r] 1 & 0 \cr 0 & 1 \end{bmatrix*}
  + \xi_0 \begin{bmatrix*}[r] 0 & \shortm 1 \cr 0 & 0 \end{bmatrix*}\right) \*x_2                                                       & = \vt{0} \label{eq:4par_main_c}\\
  \left(\begin{bmatrix*}[r] 0 & 0 \cr 0 & 1 \end{bmatrix*}
  + \I k \begin{bmatrix*}[r] 0 & 1 \cr 1 & 0 \end{bmatrix*}
  + \xi_0 \begin{bmatrix*}[r] 1 & 0 \cr 0 & 0 \end{bmatrix*}\right) \*x_3                                                       & = \vt{0} \label{eq:4par_main_d}
\end{align}
\end{subequations}
where
\begin{equation}\label{eq:xi1_xi2_xi3_def}
  \kappa_{y,1}:=\sqrt{\kappa_1^2 - k^2}, \quad \kappa_{y,2}:=\sqrt{\kappa_2^2 - k^2},\quad \xi_0=-k^2
\end{equation}
and $\*x_1\neq 0$, $\*x_2\neq 0$ and $\*x_3\neq 0$. 
\end{theorem}
\begin{proof}
The first equation \eqref{eq:4par_main_a} is satisfied by the  definition of $\kappa_{y,1}$, $\kappa_{y,2}$ and $\xi_0$ in \eqref{eq:xi1_xi2_xi3_def}. 
The determinant of the matrix in  \eqref{eq:4par_main_b}  is $\kappa_1^2+\xi_0-\kappa_{y,1}^2$, which vanishes for $\kappa_{y,1}$ and $\xi_0$ given in the definition \eqref{eq:xi1_xi2_xi3_def}. Since the matrix is singular, there exists a non-zero singular vector $\*x_1$ such that \eqref{eq:4par_main_b} is satisfied. Equation \eqref{eq:4par_main_c} follows analogously. The determinant of the matrix in \eqref{eq:4par_main_d} is $\xi_0+k^2$, which is zero for the definition of $\xi_0$, and there exists a non-zero singular vector such that \eqref{eq:4par_main_d} is satisfied.
\end{proof}

The above construction is by no means unique, i.e., there are other ways to build a multiparameter eigenvalue problem containing the solutions of the original problem. However, it is the formulation with the smallest number of additional auxiliary variables that we could find for the general case. 
Moreover, we built it such that if all matrices of the original problem (finite-element matrices and coupling matrices) are real, then 
all matrices of the multiparameter eigenvalue problem \eqref{eq:4par} are real. This is important for numerical
computations as it is well known that eigenvalue problems with real matrices can be generally solved faster than problems with complex matrices (irrespective of whether the eigenvalues are real-valued).

For special cases, there are transformations with fewer variables. In the case where the plate is either coupled to a fluid only at one of its surfaces or to the same fluid at both surfaces, we may simply eliminate one parameter, say $\kappa_{y,2}$, and remove the corresponding third equation in the above system, yielding a three-parameter eigenvalue problem. Specifically, in the case of the same fluid on both sides, we have $\kappa_{y,1} = \kappa_{y,2}$ and hence $\kappa_{y,1}\*R_1 + \kappa_{y,2} \*R_2 \eqqcolon \kappa_{y_1} \*R$. However, as mentioned before, there are more efficient ways to solve this special case~\cite{Kiefer2019}. Another important special case is the one where the plate consists only of isotropic materials, in which case the problem can be further simplified as shown in \ref{app:isotropic}.

In the case of coupling to solid media with $\*R(k)$ is given by \eqref{eq:nonpoly2}, we introduce
new parameters
\begin{equation}
  \xi_1=k \kappa_{y,1},\quad
  \xi_2=k \gamma_{y,1}, \quad
  \xi_3=k \kappa_{y,2}, \quad
  \xi_4=k \gamma_{y,2}
\end{equation}
and, together with $\I k$ and $\xi_0$, obtain the six-parameter eigenvalue problem
\begin{subequations}\label{eq:6par}
\begin{align}
  \left(-\Etbar+\omega^2\,\Mbar + \I k \Eobar+\xi_1 \*R_{1,1} + \xi_2 \*R_{2,1}+\xi_3 \*R_{1,2} + \xi_4 \*R_{2,2} + \xi_0 \Ezbar \right)\egv & =\vt{0}\label{eq:6par_main}\\
  \left(\begin{bmatrix*}[c] 0 & \shortm \kappa_{1}^2 \cr 0 & 0 \end{bmatrix*}
  + \xi_1 \begin{bmatrix*}[r] \minusspace 1 & 0 \cr 0 & \minusspace 1 \end{bmatrix*}
  + \xi_0 \begin{bmatrix*}[r] 0 & \shortm 1 \cr 1 & 0 \end{bmatrix*}
  \right) \*x_1                                                                                                                                   & = \vt{0} \label{eq:6par_main_b}\\
  \left(\begin{bmatrix*}[c] 0 & \shortm \gamma_{1}^2 \cr 0 & 0 \end{bmatrix*}
  + \xi_2 \begin{bmatrix*}[r] \minusspace 1 & 0 \cr 0 & \minusspace 1 \end{bmatrix*}
  + \xi_0 \begin{bmatrix*}[r] 0 & \shortm 1 \cr 1 & 0 \end{bmatrix*}
  \right) \*x_2                                                                                                                                   & = \vt{0} \label{eq:6par_main_c}\\
  \left(\begin{bmatrix*}[c] 0 & \shortm \kappa_{2}^2 \cr 0 & 0 \end{bmatrix*}
  + \xi_3 \begin{bmatrix*}[r] \minusspace 1 & 0 \cr 0 & \minusspace 1 \end{bmatrix*}
  + \xi_0 \begin{bmatrix*}[r] 0 & \shortm 1 \cr 1 & 0 \end{bmatrix*}
  \right) \*x_3                                                                                                                                   & = \vt{0} \label{eq:6par_main_d}\\
  \left(\begin{bmatrix*}[c] 0 & \shortm \gamma_{2}^2 \cr 0 & 0 \end{bmatrix*}
  + \xi_4 \begin{bmatrix*}[r] \minusspace 1 & 0 \cr 0 & \minusspace 1 \end{bmatrix*}
  + \xi_0 \begin{bmatrix*}[r] 0 & \shortm 1 \cr 1 & 0 \end{bmatrix*}
  \right) \*x_4                                                                                                                                   & = \vt{0} \label{eq:6par_main_e}\\
  \left(\begin{bmatrix*}[r] 0 & 0 \cr 0 & 1 \end{bmatrix*}
  + \I k \begin{bmatrix*}[r] 0 & 1 \cr 1 & 0 \end{bmatrix*}
  + \xi_0 \begin{bmatrix*}[r] \minusspace 1 & 0 \cr 0 & 0 \end{bmatrix*}\right) \*x_5                                                                                                 & = \vt{0} \label{eq:6par_main_f}
\end{align}
\end{subequations}
with the definitions given in Eq.~\eqref{eq:xi1_xi2_xi3_def} and, additionally,
\begin{equation}\label{eq:gammay_def}
  \gamma_{y,1}:=\sqrt{\gamma_1^2 - k^2}, \quad \gamma_{y,2}:=\sqrt{\gamma_2^2 - k^2}.
\end{equation}
Note that the determinant of the matrix in Eq.~\eqref{eq:6par_main_b} (analogously for \eqref{eq:6par_main_c}--\eqref{eq:6par_main_e}) is now $\xi_0(\kappa_1^2+\xi_0-\kappa_{y,1}^2)$, which again vanishes for $\kappa_{y,1}$ and $\xi_0$ given by \eqref{eq:xi1_xi2_xi3_def}, but also for $k=0$ (rigid body modes).
Again, if the plate is coupled to the same material on both sides or one of the surfaces is free, we can reduce the number of parameters. In this case, we can remove $\xi_3$ and $\xi_4$ and the corresponding fourth and fifth equations, resulting in a four-parameter eigenvalue problem. We can also account for the situation where the plate is coupled to a fluid medium at one surface and a solid at the other. In this case, we include one parameter for the fluid, say $\xi_1$, as in Eq.~\eqref{eq:4par} and two parameters, say $\xi_3$, $\xi_4$, in the form of \eqref{eq:6par}, resulting in a five-parameter eigenvalue problem.

The multiparameter eigenvalue problems \eqref{eq:4par} and \eqref{eq:6par} are both singular. It is known that we can find eigenvalues using a generalized staircase-type numerical algorithm from \cite{MP_Q2EP}, for which a Matlab implementation is available \cite{multipareig_2023}. Hence, we will not explain the numerical computation in much detail. In \ref{appendix}, we present a variant that works for the 
particular problems \eqref{eq:4par} and \eqref{eq:6par} and uses a shift to avoid additional steps necessary to deal with singular problems, resulting in a more efficient computation.

In general, $\I k$ is an eigenvalue of a generalized eigenvalue problem 
$\vg{\Delta}_1 \*v=\I k \vg{\Delta}_0 \*v$, where matrices
$\vg{\Delta}_0$ and $\vg{\Delta}_1$ are the so-called operator determinants. These matrices are of size $n_\Delta \times n_\Delta$ with $n_\Delta = \prod \limits_{i=1}^r n_i$, i.e., the product of the sizes of the matrices in each equation, cf., \eqref{eq:mupa}.
Specifically, the operator determinants are of size $8n_1\times 8n_1$ for \eqref{eq:4par} and $32n_1\times 32n_1$ for \eqref{eq:6par}, where $n_1\times n_1$ is the size of the matrices in \eqref{eq:4par_main_a} or \eqref{eq:6par_main}. Hence, the computational effort increases rapidly with the number of sought parameters.
The details on the construction of $\vg{\Delta}_0$ and $\vg{\Delta}_1$ 
from the matrices in \eqref{eq:4par} and \eqref{eq:6par} are summarized in \ref{sec:appx_mep}, see also \cite{Kiefer2022} or \cite{PCH_SpecColl_JCP}. Note that, by solving the multiparameter eigenvalue problem \eqref{eq:4par}, we obtain not only the eigenvalues $k$ but also the values of the corresponding parameters $\kappa_{y,1}$ and $\kappa_{y,2}$ and, thus, the correct signs of $\sqrt{\kappa_1^2-k^2}$ and $\sqrt{\kappa_2^2-k^2}$. The same applies to \eqref{eq:6par}, where we get, in addition, the values of $\sqrt{\gamma_1^2-k^2}$ and $\sqrt{\gamma_2^2-k^2}$.

A few remarks on the above formulation are in order. Firstly, in the notation adopted here, $k^2$ is identified as an additional parameter. Alternatively, we could linearize Eqs.~\eqref{eq:4par} and \eqref{eq:6par_main} in $k$ with a companion linearization, which doubles the size of the corresponding matrices. This is usually done when solving the quadratic eigenvalue problem corresponding to the free plate (Section~\ref{sec:safe}). Both variants lead to operator determinants of the same size. However, the first variant offers the advantage that it can be written in such a way that all matrices in the multiparameter eigenvalue problems are real (provided that the finite-element matrices are real), making the computation more efficient.
In the case of coupling to fluid media, the parameters directly correspond to $k$, $k^2$, $\kappa_{y,1}$, $\kappa_{y,2}$, the last two being the vertical wavenumbers in the two fluids. In contrast, coupling to solid media results in additional parameters of the form $k\kappa_y$ and $k\gamma_y$. 

We solve the multiparameter eigenvalue problem for a set of given frequencies; hence, the frequency $\omega$ is treated as a constant rather than an additional parameter.
Furthermore, it is important to note that the solution procedure allows both positive and negative signs of all wavenumbers, including those of any (partial) waves in the half-spaces. That is to say, we do not only obtain solutions characterized by waves propagating \textit{away} from the plate inside the unbounded media, but also \textit{towards} it. For a given physical setup, the relevant modes need to be selected based on the combination of signs in a postprocessing step. Our numerical studies focus on those scenarios where guided waves propagate inside a plate structure and radiate energy into the adjacent half-spaces. However, it can be interesting to study the other types of modes in the context of, e.g., air-coupled or immersed ultrasonic testing or seismic waves interacting with soil layers.

\section{Numerical examples}\label{sec:examples}\noindent
In order to validate and illustrate the competitiveness of our approach, we present four numerical examples that involve plate structures of different materials coupled to acoustic and elastic half-spaces. Homogenous elastic plates are discretized by a single finite element whose polynomial degree $p_e$ is chosen according to the recommendation given in \cite{Gravenkamp2014} as
\begin{equation}\label{eq:orderChoice}
  p_e = \frac{a_0}{2} + 3, \qquad  a_0 = \frac{h\,\omega_\mathrm{max}}{c_s}.
\end{equation}
Here, $a_0$ denotes a dimensionless frequency, computed based on the shear wave velocity $c_s$, layer thickness $h$, and maximum frequency of interest $\omega_\mathrm{max}$. When modeling plates consisting of several layers, the polynomial degree of each layer is chosen separately based on its local shear wave velocity and thickness. For each example, we plot the dispersion curves in terms of phase velocity and attenuation. In addition, we visualize the wave fields created by two selected modes at different frequencies in order to provide some physical insight into the wave propagation behavior. The relevant parameters of all materials used in the examples are summarized in Table~\ref{tab:materialParameters}. All computations have been performed using a Matlab implementation of the proposed approach, executed on a desktop computer using an 11$^\text{th}$ Gen Intel Core i9-11900K processor (3.50 GHz) and 16 GB RAM.

\begin{table} \centering
  \caption{Overview of material parameters used in the numerical experiments. \label{tab:materialParameters}}
  \begin{tabular*}{0.7\textwidth}[tb]{c|cccrr}
            & density $\rho$ & \multicolumn{2}{c}{ wave speeds $c_\ell, c_t$ } &  \multicolumn{2}{c}{ Lam\'e parameters $\lambda, \mu$}\\ \hline
    brass   & 8.40 g/cm$^3$ & 4.40 km/s & 2.20 km/s & 81.312 GPa &  40.656 GPa \\
    Teflon  & 2.20 g/cm$^3$ & 1.35 km/s & 0.55 km/s & 2.679 GPa & 0.666 GPa \\
    titanium & 4.46 g/cm$^3$ & 6.06 km/s & 3.23 km/s & 70.726 GPa & 46.531 GPa \\
    water   & 1.00 g/cm$^3$ & 1.48 km/s &         &  \\
    oil     & 0.87 g/cm$^3$ & 1.74 km/s &         &
  \end{tabular*}
\end{table}

\subsection{Brass plate immersed in water}
\begin{figure}\centering
  \subfloat[]{\includegraphics[width=0.5\textwidth]{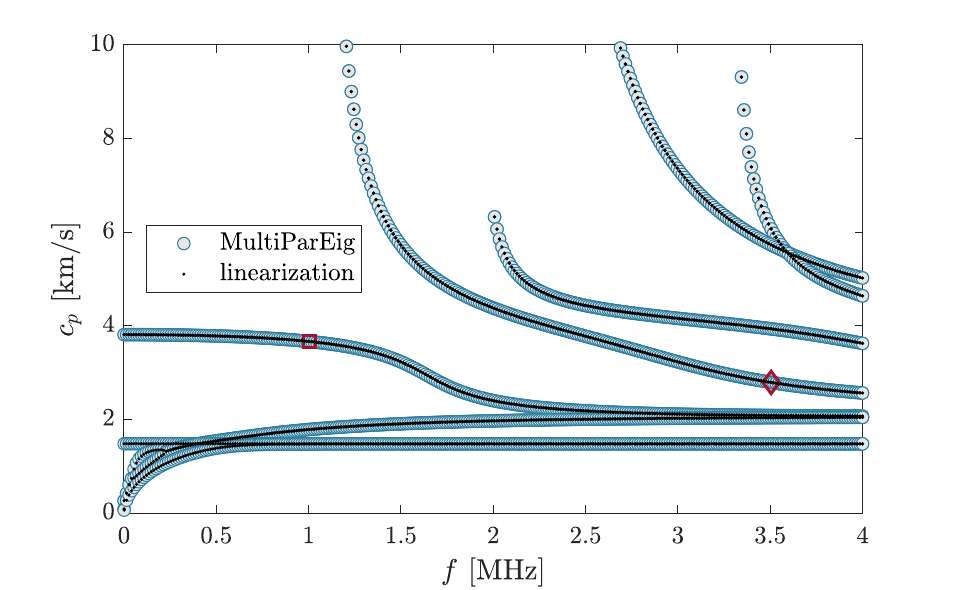}}
  \subfloat[]{\includegraphics[width=0.5\textwidth]{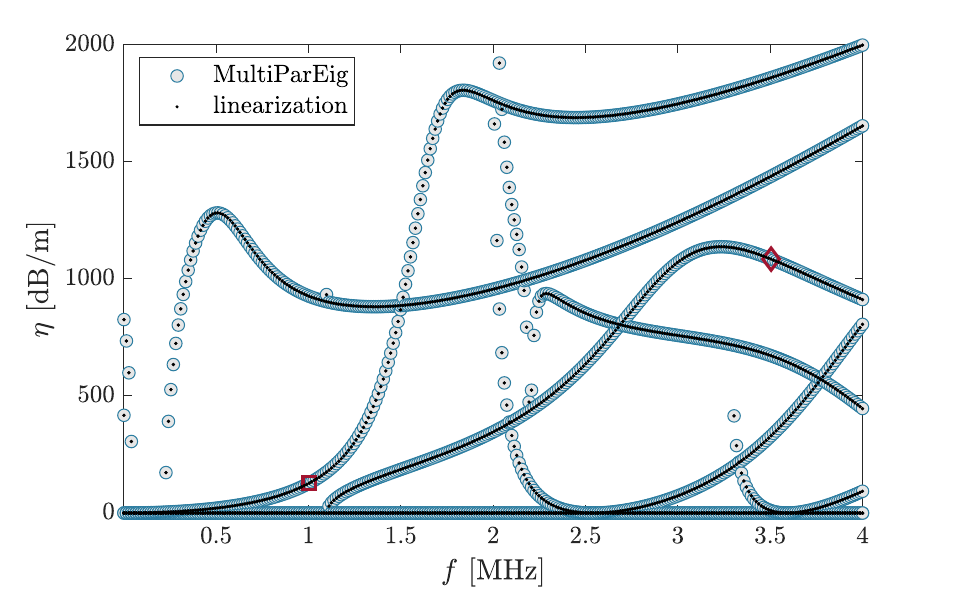}}
  \caption{(a) Phase velocity and (b) attenuation of guided waves in a 1\,mm thick brass plate immersed in water. Results are computed using the proposed method (`MultiParEig') and validated against the alternative approach based on a linearization~\cite{Kiefer2019}. The symbols \textcolor{myred}{$\bm{\square}$} and \textcolor{myred}{$\bm{\Diamond}$} indicate arbitrarily chosen modes whose wave fields are plotted in Fig.~\ref{fig:brassWater_modeshapes} \label{fig:brassWater_curves} }
\end{figure}

\noindent
As a first numerical experiment for validating the proposed approach, we consider waves propagating along a brass plate of thickness 1\,mm immersed in water. We compute dispersion properties for 300 frequencies up to 4\,MHz. According to Eq.~\eqref{eq:orderChoice} and considering the material parameters presented in Table~\ref{tab:materialParameters}, a discretization using one element with a polynomial degree of 9 (10 nodes) is sufficient to obtain accurate results. We assume plane strain conditions, hence considering in-plane displacements only (`Lamb-type' modes). The shear-horizontal (out-of-plane) modes are omitted here as they do not couple to the acoustic medium.  
Consequently, we obtain finite-element matrices of size $22\times22$, including the two DOFs representing the pressure at the upper and lower surface.
As the plate is coupled to the same fluid on both sides, we can employ, for validation, the simpler approach of linearizing the eigenvalue problem by means of a change of variables as described in \cite{Kiefer2019}. Figure~\ref{fig:brassWater_curves} shows the dispersion curves in terms of phase velocity and attenuation. The approach proposed in this paper is labeled `MultiParEig' (multiparameter eigenvalue problem), while `linearization' refers to the method in \cite{Kiefer2019}. Results are in excellent agreement between both approaches. 
The computing time required for obtaining the complete dispersion curves was about 1.6\,s (i.e., 5\,ms per frequency) when exploting the linearization described in \ref{app:isotropic} and 3.6\,s without this improvement. 

\begin{figure}\centering
  \subfloat[ \label{fig:brassWater_modeshapes_smallF}]{\includegraphics[width=0.49\textwidth]{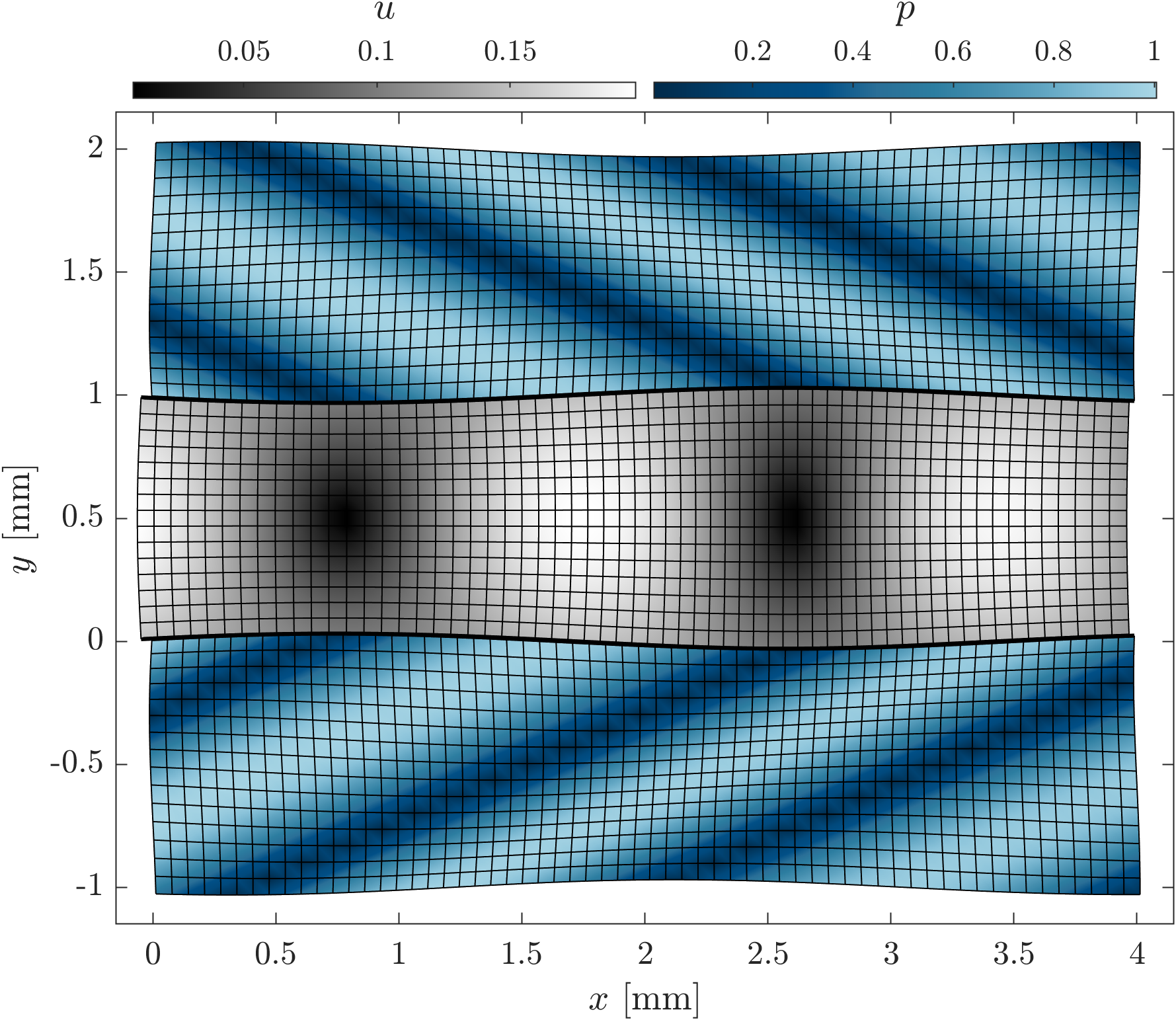}}\hfill
  \subfloat[ \label{fig:brassWater_modeshapes_largeF}]{\includegraphics[width=0.49\textwidth]{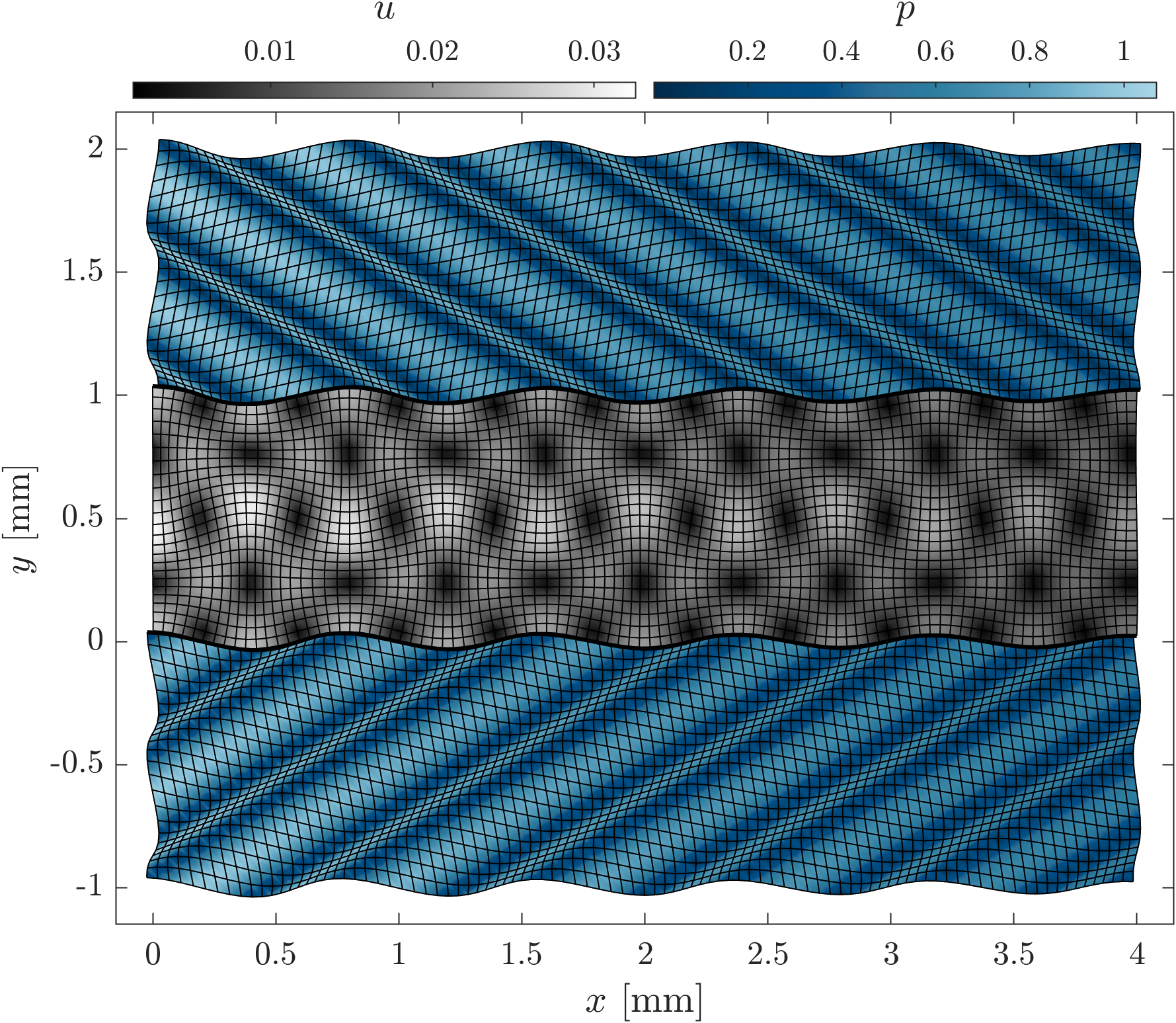}}
  \caption{Wave fields of two modes in the water-immersed brass plate at (a) 1\,MHz and (b) 3.5\,MHz. The modes are indicated in the dispersion diagram in Fig.~\ref{fig:brassWater_curves} by the symbols \textcolor{myred}{$\bm{\square}$} and \textcolor{myred}{$\bm{\Diamond}$}.    \label{fig:brassWater_modeshapes}  }
\end{figure}

In Fig.~\ref{fig:brassWater_modeshapes}, we present the wave fields of two arbitrarily selected modes at 1\,MHz and 3.5\,MHz. The modes are indicated in the dispersion graphs as {\small\textcolor{myred}{$\pmb{\Box}$}} and {\small\textcolor{myred}{$\pmb{\Diamond}$}}. Note that, even though we only discretize the plate's thickness, the wave field anywhere inside the plate as well as in the surrounding medium can be evaluated in a post-processing step. In particular, the displacements inside the plate for arbitrary values of $x$ are obtained by Eq.~\eqref{eq:plane_wave_ansatz} and interpolation along the thickness using the finite-element trial functions. Similarly, the acoustic pressure in the water is evaluated based on Eqs.~\eqref{eq:acousticPlaneWave} with the vertical wavenumbers obtained as solutions to the multiparameter eigenvalue problem.
Furthermore, the particle displacement in the acoustic domains is obtained as 
\begin{equation}
  \*u_p = \tfrac{1}{\omega^2\,\rho}\, \nabla p.
\end{equation}
In Fig.~\ref{fig:brassWater_modeshapes}, the acoustic particle displacement is utilized to show the grid distortion, while the color scale represents the acoustic pressure. In the plate, the color scale indicates the local displacement magnitude. As expected, it can be observed that the horizontal wavenumber in the acoustic domains matches that of the guided wave inside the plate. The horizontal displacement is discontinuous at the material interfaces since the transversal wave components are not transmitted into the fluid. At the higher frequency where the wavelength is small and attenuation is large (Fig.~\ref{fig:brassWater_modeshapes_largeF}), we can also note how the amplitude decreases as the wave propagates in positive $x$-direction.

\subsection{Brass plate coupled to infinite Teflon}
\begin{figure}\centering
  \subfloat[]{\includegraphics[width=0.5\textwidth]{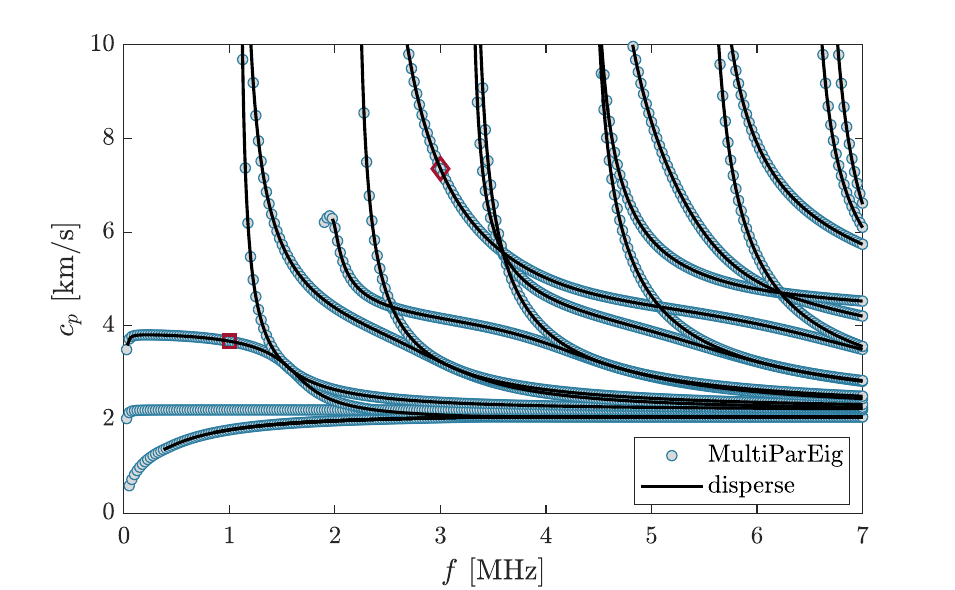}}
  \subfloat[]{\includegraphics[width=0.5\textwidth]{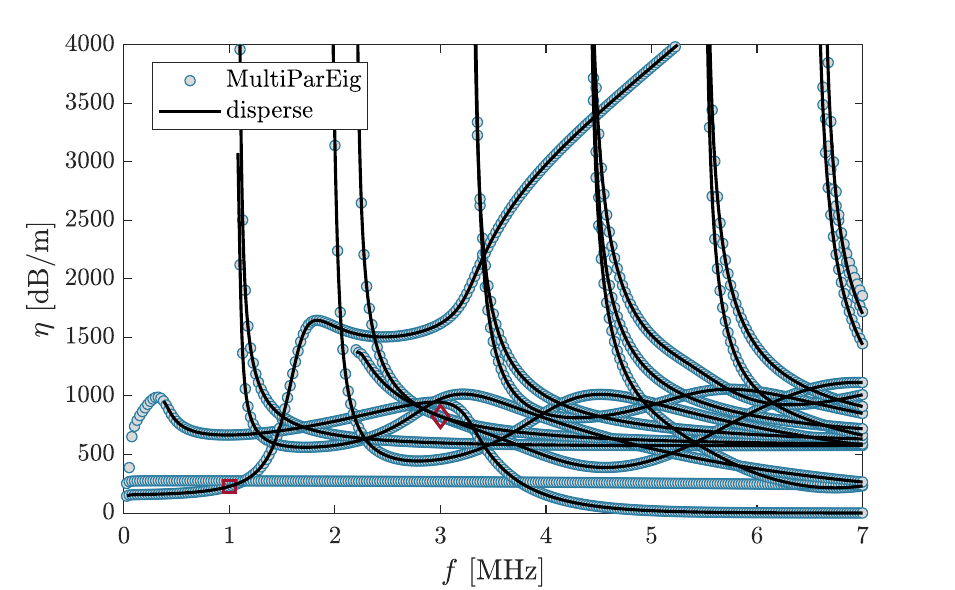}}
  \caption{(a) Phase velocity and (b) attenuation of guided waves in a 1-mm-thick brass plate coupled on one side to infinite Teflon. Results are computed using the proposed FEM-based approach and validated against the software \textit{disperse}. \label{fig:brassTeflon_curves} }
\end{figure}

\begin{figure}\centering
  \subfloat[1\,MHz]{\includegraphics[width=0.48\textwidth]{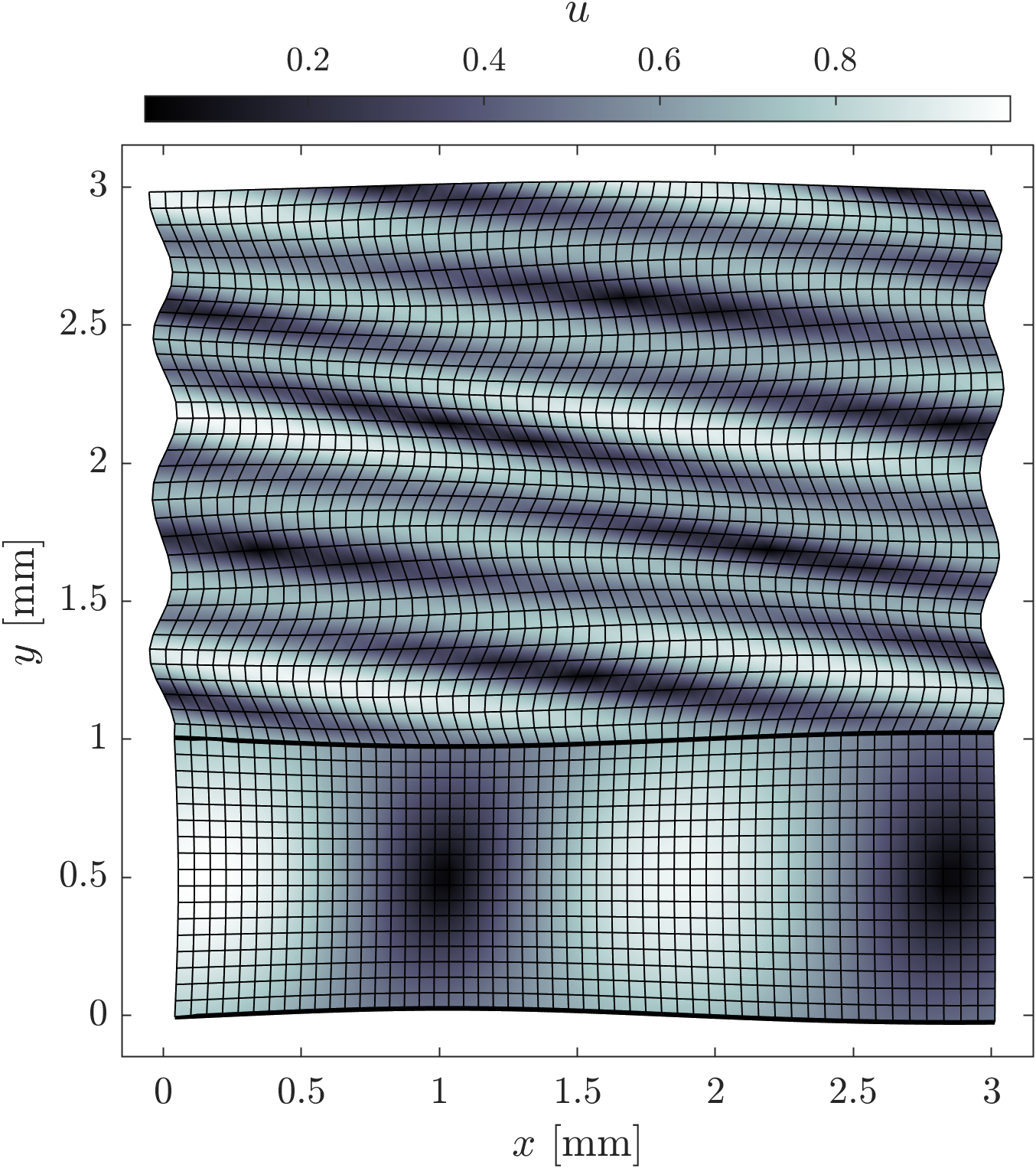}}\hfill
  \subfloat[3\,MHz]{\includegraphics[width=0.48\textwidth]{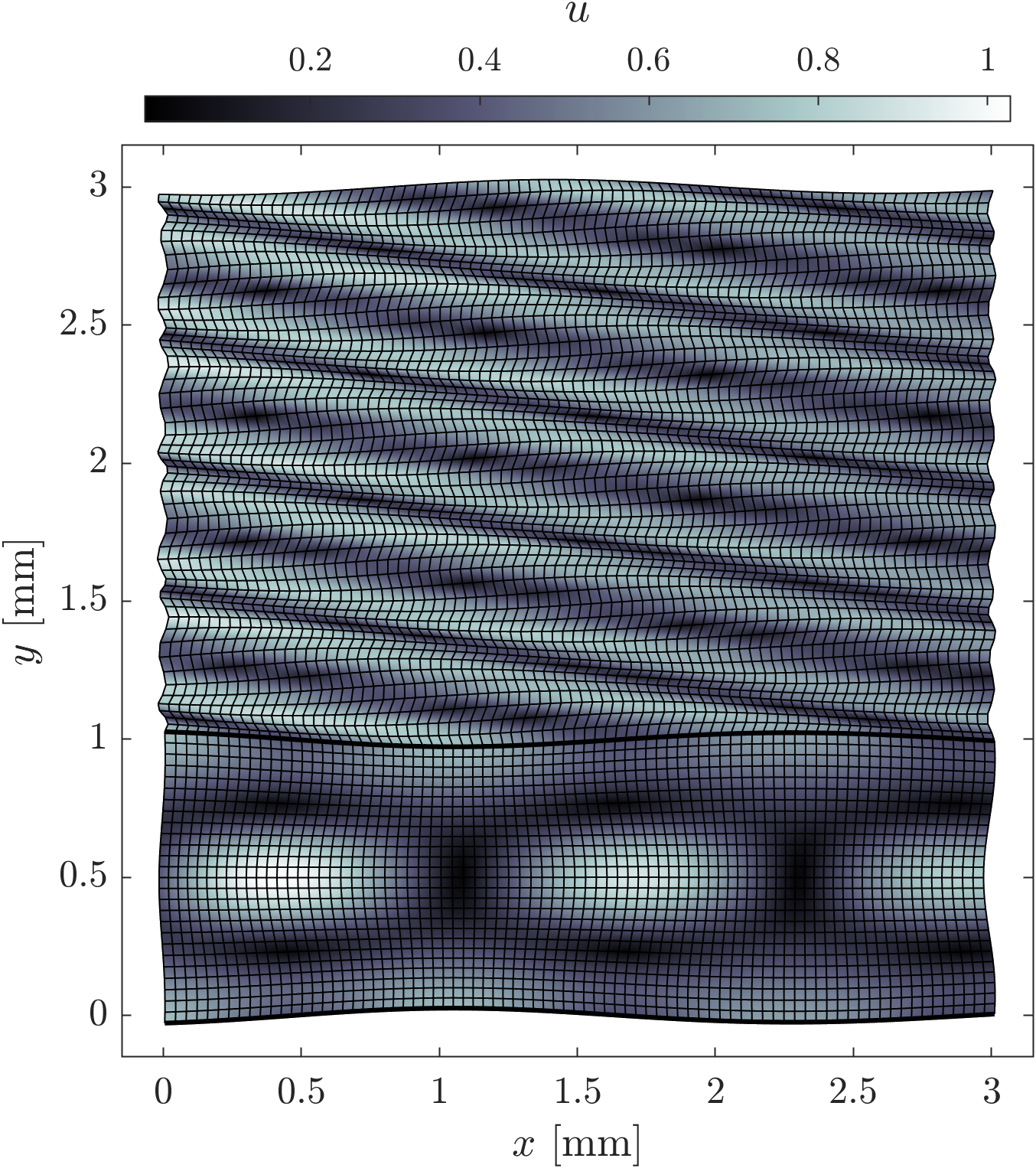}}
  \caption{Mode shapes in a brass plate coupled to Teflon of (a) the fundamental symmetric mode at 1\,MHz and (b) the seventh propagating mode at 3\,MHz. The modes are indicated in the dispersion diagram in Fig.~\ref{fig:brassTeflon_curves} by the symbols \textcolor{myred}{$\bm{\square}$} and \textcolor{myred}{$\bm{\Diamond}$}.  \label{fig:brassTeflon_modeshapes} }
\end{figure}

\noindent
To validate the computation of waveguide modes radiating into another elastic medium, we consider the same 1-mm-thick brass plate as in the previous example, now attached on one surface to a half-space with the elastic parameters of Teflon, as given in Table~\ref{tab:materialParameters}. 
The choice of this rather arbitrary academic example is motivated by the circumstance that it can be verified using the free demonstration version of the commercial software \textit{disperse} \cite{Pavlakovic2011, Pavlakovic1997}. 
Here, we compute dispersion curves up to a frequency of 7\,MHz in increments of 25\,kHz; hence, we require a finite element of order 13. In this example, we include not only the in-plane (`Lamb-modes'), but also the shear-horizontal modes that cause displacements in the out-of-plane direction $z$ only. Consequently, we obtain matrices of size $45\times45$, including three DOFs for the unbounded domain.
Phase velocity and attenuation of all leaky modes are shown in Fig.~\ref{fig:brassTeflon_curves}. The results are in excellent agreement between the proposed approach and the Global Matrix Method implemented in the disperse software. We obtained a few additional solutions where the automatic tracing algorithm in disperse failed to find the roots; apart from that, there is no discernible difference between the methods.  The computational costs for this example are significantly larger compared to the simpler case in the previous subsection due to the larger number of parameters and larger matrix sizes, requiring about 35\,s to compute the shown dispersion curves (123\,ms per frequency).
Again, we visualize in Fig.~\ref{fig:brassTeflon_modeshapes} the mode shapes of two modes at frequencies of 1\,MHz and 3.5\,MHz by computing the displacement field in a section of the plate and part of the adjacent medium. Since both materials are now elastic, both longitudinal and shear components are transmitted from the plate into the unbounded medium, and consequently, all displacement components are continuous at the interface.

\subsection{Solid layer between two elastic half-spaces}
\begin{figure}\centering
  \subfloat[]{\includegraphics[width=0.5\textwidth]{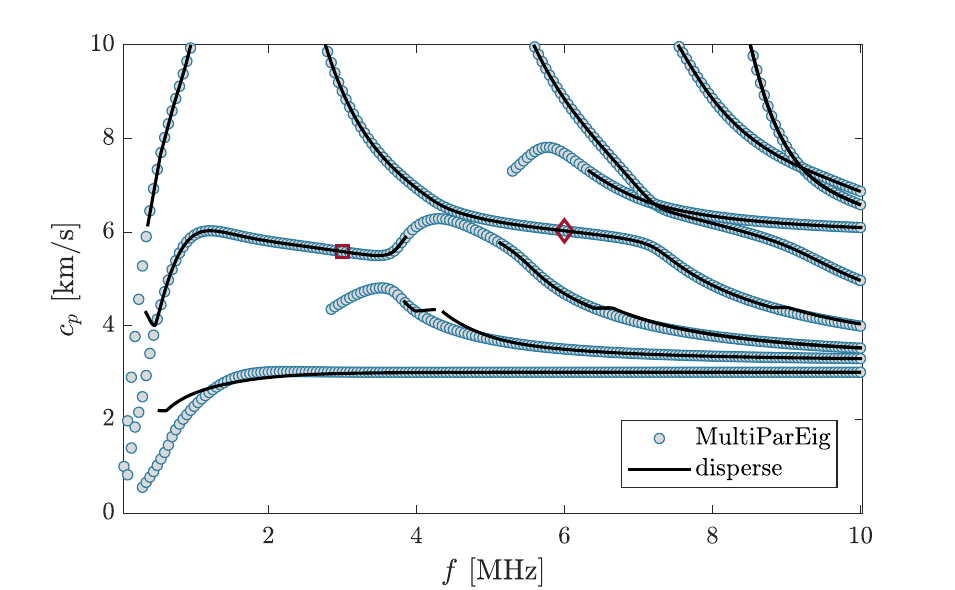}}
  \subfloat[]{\includegraphics[width=0.5\textwidth]{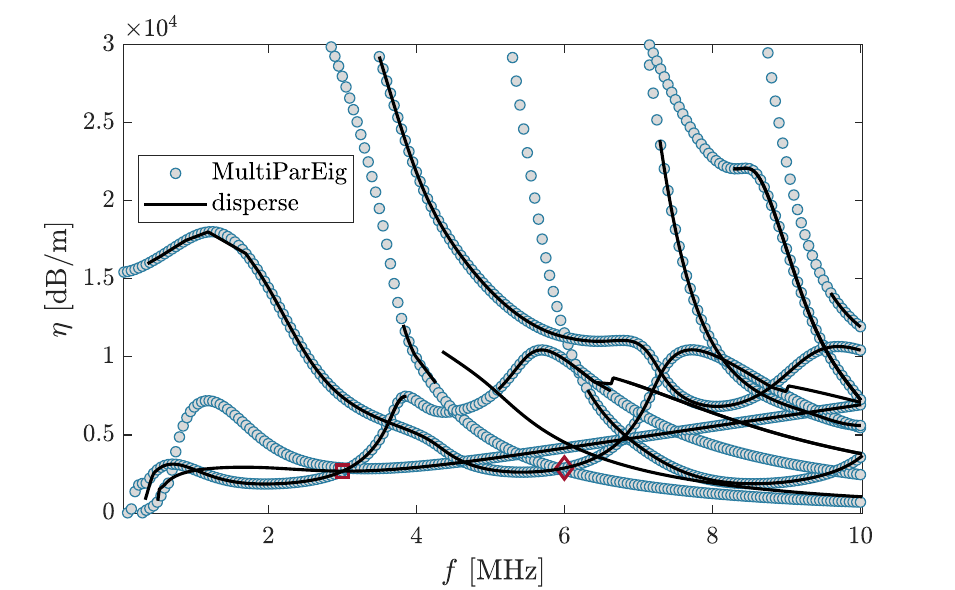}}
  \caption{(a) Phase velocity and (b) attenuation of guided wave modes in a 1\,mm thick titanium plate between two elastic half-spaces consisting of brass and Teflon, respectively. Results are computed using the proposed FEM-based approach and validated against the software \textit{disperse}. \label{fig:titanium_brassTeflon_curves} }
\end{figure}

\begin{figure}\centering
  \subfloat[3\,MHz]{\includegraphics[width=0.48\textwidth]{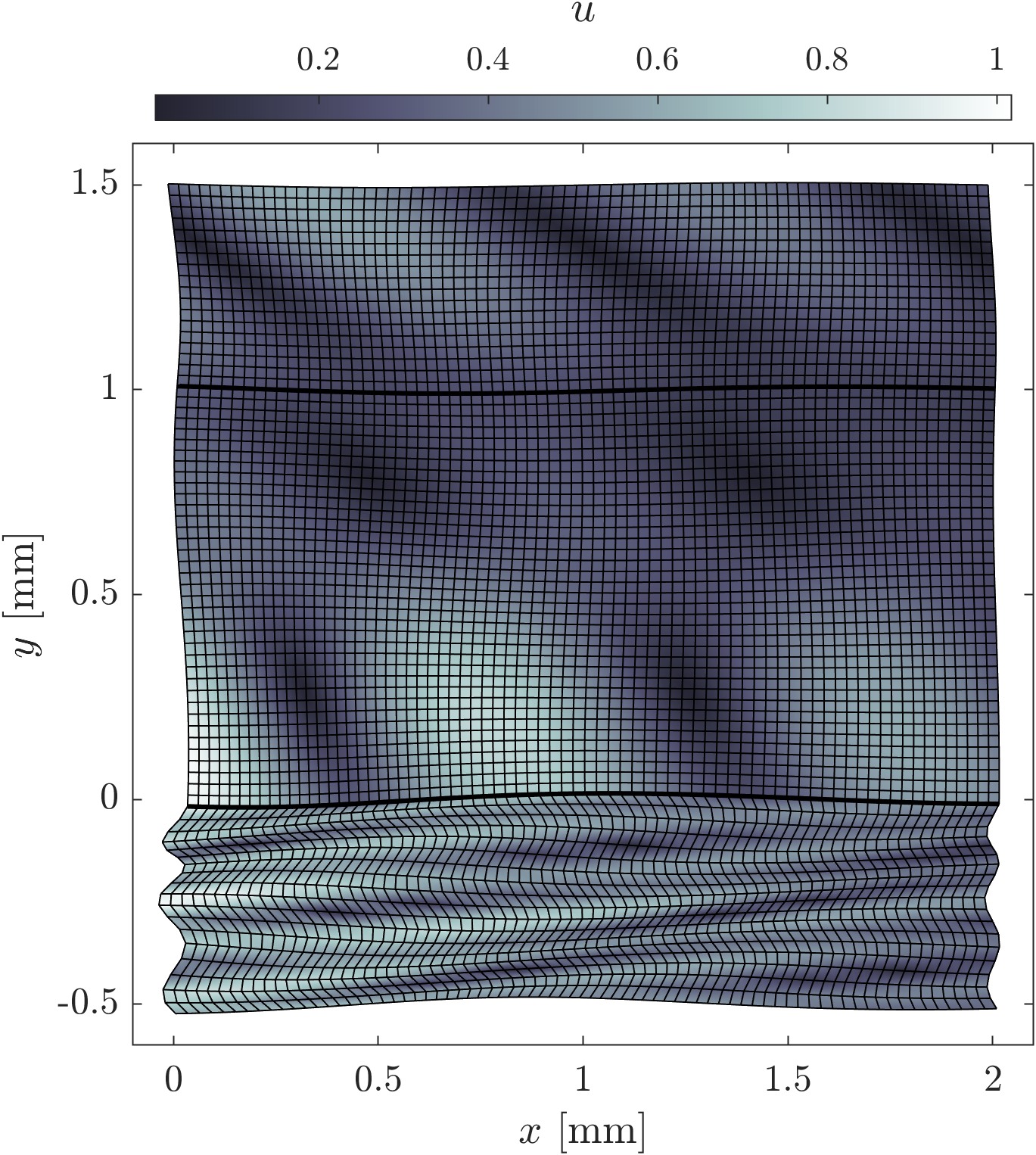}}\hfill
  \subfloat[6\,MHz]{\includegraphics[width=0.48\textwidth]{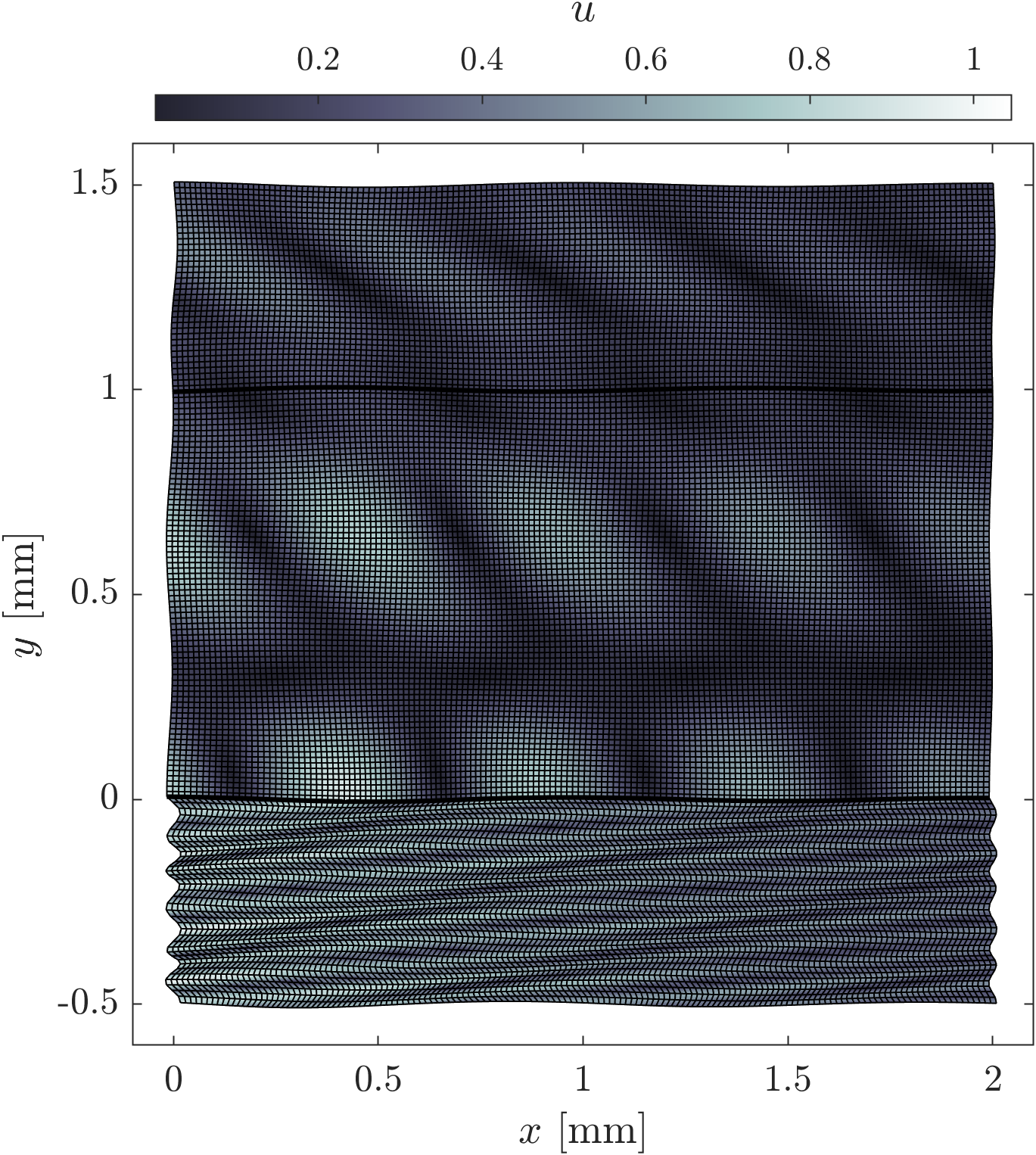}}
  \caption{Mode shapes of two modes at (a) 3\,MHz and (b) 6\,MHz for the titanium plate between a Teflon half-space (below) and a brass half-space (above). The modes are indicated in the dispersion diagram in Fig.~\ref{fig:titanium_brassTeflon_curves} by the symbols \textcolor{myred}{$\bm{\square}$} and \textcolor{myred}{$\bm{\Diamond}$}.    \label{fig:titanium_brassTeflon_modeshapes} }
\end{figure}

\noindent
In this example, we address the computationally most demanding case in which the plate is coupled to two different solid media at the lower and upper surfaces, leading to the six-parameter eigenvalue problem in Eq.~\eqref{eq:6par}. The layer consists of titanium, while the elastic parameters of Teflon (below the plate) and brass (above the plate) are assumed for the half-spaces. We choose a frequency range up to 10\,MHz, requiring a finite element of order 13 to discretize the plate. The frequency increment was selected as 50\,kHz. 
Results for the `Lamb-type' modes are presented in Fig.~\ref{fig:titanium_brassTeflon_curves}; again, a comparison is made with the software \textit{disperse}. This example is particularly challenging due to the relatively small acoustic mismatch between titanium and brass, resulting in significant radiation into the upper half-space. When using disperse, -- which requires tracing individual modes starting from isolated roots of the dispersion relation -- finding all solutions is not straightforward. To achieve the shown comparison, we made extensive use of the software's feature to search for roots in regions where we expected to find modes based on our own computations. Still, we were not able to trace all branches over the entire frequency range in disperse, and for some modes, the tracing algorithm failed to find the correct imaginary part of the solution. Nevertheless, the modes that were traced successfully agree well with our solutions, which shall suffice to demonstrate the validity of both formulations. Our approach required about 186\,s to compute the dispersion curves (920\,ms per frequency).
The wave fields visualized in Fig.~\ref{fig:titanium_brassTeflon_modeshapes} help demonstrate the difficulty of this example. Due to the drastically different elastic constants between titanium and teflon, waves propagate with significantly smaller vertical wavelengths in the lower half-space. At the same time, the difference between titanium and brass is so small that the interface at $y=1\,\text{mm}$ is hardly noticeable. However, even for this rather extreme case, the proposed approach yields robust results.

\subsection{Layered plate between a fluid and an elastic half-space}
\begin{figure}\centering
  \subfloat[]{\includegraphics[width=0.49\textwidth]{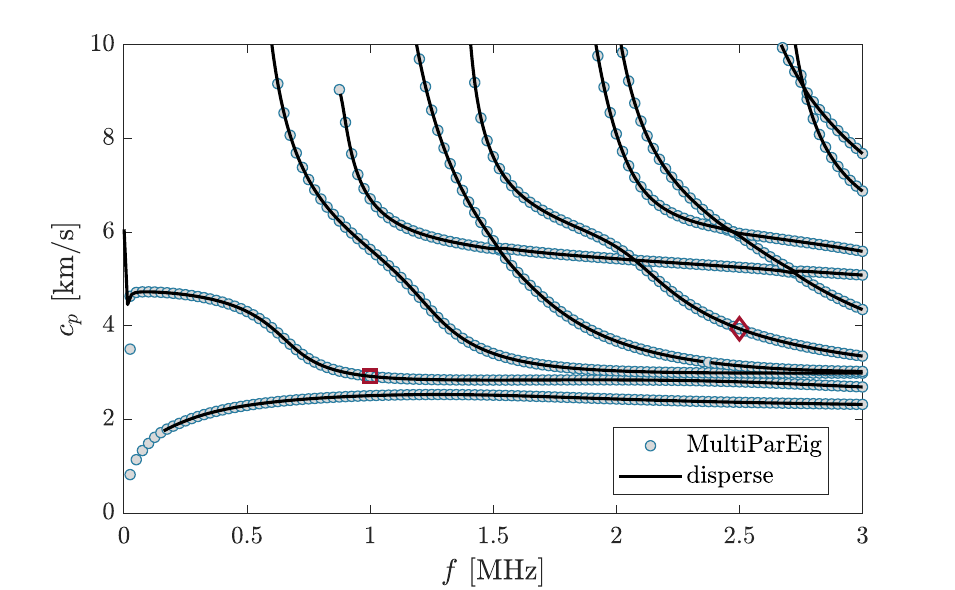}}\hfill
  \subfloat[]{\includegraphics[width=0.49\textwidth]{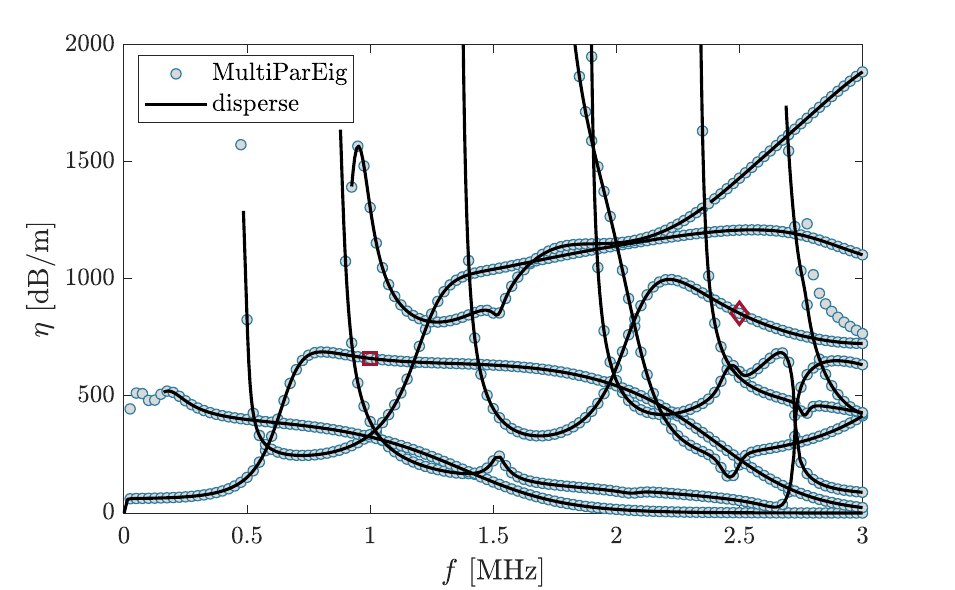}}
  \caption{(a) Phase velocity and (b) attenuation of guided wave modes in a layered plate (titanium-brass-titanium) between an elastic (Teflon) and a fluid (oil) half-space. Results are computed using the proposed FEM-based approach and validated against the software \textit{disperse}. \label{fig:composite_oilTeflon_curves} }
\end{figure}

\begin{figure}\centering
  \subfloat[1\,MHz]{\includegraphics[width=0.48\textwidth]{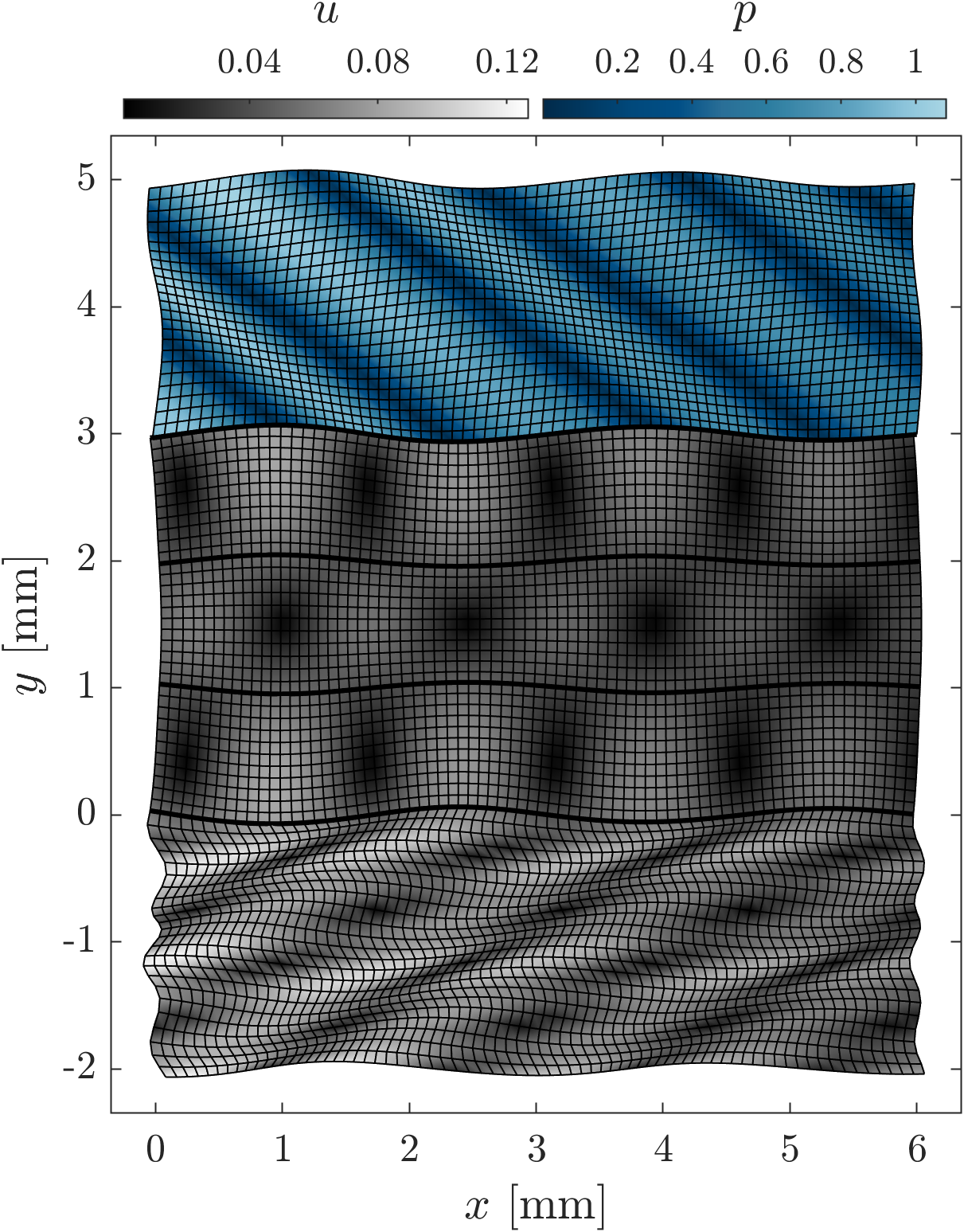}}\hfill
  \subfloat[2.5\,MHz]{\includegraphics[width=0.48\textwidth]{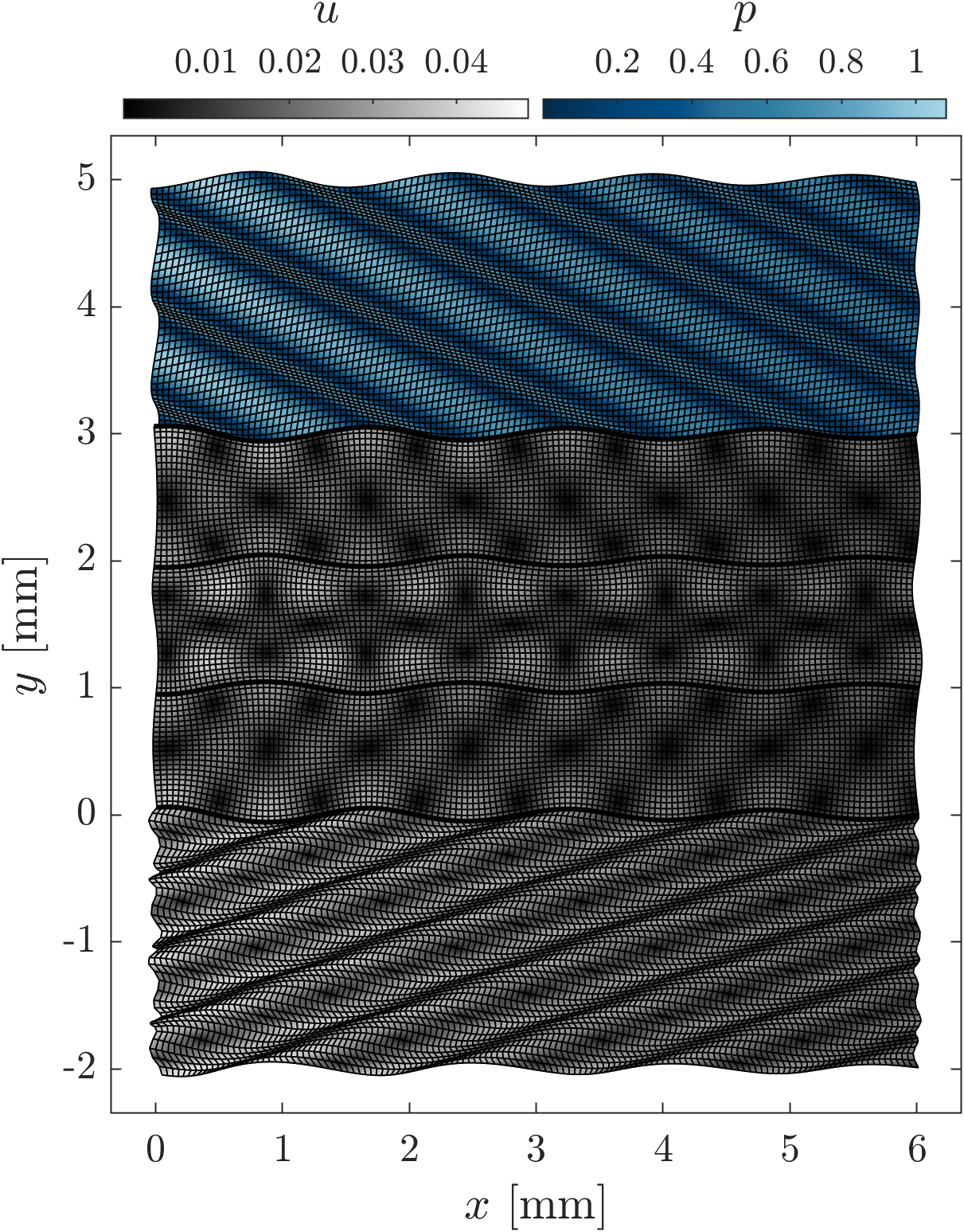}}
  \caption{Mode shapes of two waves at (a) 1\,MHz and (b) 2.5\,MHz for the trilayered plate (titanium-brass-titanium) between an elastic (Teflon) and a fluid (oil) half-space. The modes are indicated in the dispersion diagram in Fig.~\ref{fig:titanium_brassTeflon_curves} by the symbols \textcolor{myred}{$\bm{\square}$} and \textcolor{myred}{$\bm{\Diamond}$}.    \label{fig:composite_TiBrTi_unbTeflonOil_modeshapes} }
\end{figure}

\noindent
As a final numerical study, we consider a multi-layered plate coupled to an elastic and a fluid half-space at the lower and upper surfaces, respectively. The plate itself consists of three layers (titanium-brass-titanium), each with a thickness of 1\,mm. For the half-spaces, material parameters of Teflon and oil are assumed, as given in Table~\ref{tab:materialParameters}. For the brass layer, we require an element order of eight, while an order of six is sufficient for the titanium layers up to a frequency of 3\,MHz.
Assuming again a plane strain state with two displacement components per node, we obtain finite-element matrices of size $45\times 45$, including the three DOFs (two displacements and one pressure) describing the unbounded domains. 
Dispersion curves are given in Fig.~\ref{fig:composite_oilTeflon_curves}, computed at 121 frequencies with an increment of 25\,kHz. The computing time was about 56\,s (460\,ms per frequency). Again, results agree well with the reference solution obtained using \emph{disperse}, with similar challenges in the application of the Global Matrix Method as described in the previous example. Figure~\ref{fig:composite_TiBrTi_unbTeflonOil_modeshapes} illustrates the relatively complex mode shapes inside the multi-layered specimen compared to the previous examples.

\section{Conclusions}\label{sec:conclusion}\noindent
We have demonstrated that the interface conditions arising from the coupling between an elastic layered medium and adjacent acoustic or linearly elastic unbounded domains can be incorporated rigorously into a conventional semi-analytical waveguide model. The resulting nonlinear eigenvalue problem is solved by exploiting the connection to multiparameter eigenvalue problems, for which established algorithms are available. As a result, we obtain the wavenumbers of guided waves, as well as those in the adjacent media, together with the corresponding discretized mode shapes. These problems are computationally relatively demanding to solve compared to the modeling of waves in free plates. This is particularly true in the case of two adjacent elastic media, which require the solution of generalized eigenvalue problems more than 32 times larger than that describing the free plate. Nevertheless, we demonstrated that the computational costs are manageable for relevant applications when using efficient discretization schemes to model the plate. 
As of now, the presented formulation of unbounded elastic media is restricted to isotropic undamped materials. This limitation may be possible to overcome by generalizations to be done in future work.

\section*{Acknowledgments}\noindent
Hauke Gravenkamp acknowledges grant CEX2018-000797-S funded by the Ministerio de Ciencia e Innovaci\'on, MCIN/AEI/ 10.13039/501100011033.
Bor Plestenjak has been supported by the Slovenian Research and Innovation Agency (grants N1-0154 and P1-0194). Daniel A.\ Kiefer has received support under the program ``Investissements d'Avenir'' launched by the French Government under Reference No.\ ANR-10-LABX-24.

\appendix
\section{Multiparameter eigenvalue problems and operator determinants}
\label{sec:appx_mep}\noindent
An \emph{$r$-parameter eigenvalue problem} has the form
\begin{equation}\label{eq:add_mupa}
  (\*A_{i0}+\lambda_1 \*A_{i1} + \cdots + \lambda_r \*A_{ir})\,\*x_i=\vt{0},\qquad i=1,\ldots,r,
\end{equation}
where $\*A_{ij}$ is an $n_i\times n_i$ matrix and $\*x_i\ne\vt{0}$ for $i=1,\ldots,r$.
If \eqref{eq:add_mupa} is satisfied
then the tuple $(\lambda_1,\ldots,\lambda_r)$ is an \emph{eigenvalue} and 
$\*z = \*x_1\otimes \cdots \otimes \*x_r$ is the corresponding \emph{eigenvector}, where ``$\otimes$'' denotes the Kronecker product.\footnote{The Kronecker product $\*A\otimes \*B$ of the $m \times n$-matrix $\*A = [A_{ij}]$ with the $p \times q$-matrix $\*B$ yields the block matrix $[A_{ij} \*B]$ of size $m p \times n q$. For three
matrices, $\*A\otimes \*B\otimes \*C=(\*A\otimes \*B)\otimes \*C = \*A\otimes (\*B\otimes \*C)$. Vectors are treated like matrices with one column.} The problem
\eqref{eq:add_mupa}
is related to a system of $r$ generalized eigenvalue problems
\begin{equation}\label{eq:genEVP}
    \vg{\Delta}_i \*z =\lambda_i \vg{\Delta}_0 \*z,\qquad i=1,\ldots,r,
\end{equation}
where the matrices
\begin{equation}\label{eq:opdet}\vg{\Delta}_0=\left|\begin{matrix}\*A_{11} & \cdots & \*A_{1r}\cr
\vdots &  & \vdots \cr 
\*A_{r1} & \cdots & \*A_{rr}\end{matrix}\right|_\otimes
=\sum_{\sigma\in S_r}{\rm sgn}(\sigma) \, \*A_{1\sigma_1}\otimes \*A_{2\sigma_2}\otimes \cdots \otimes \*A_{r\sigma_r}
\end{equation}
and
$$\vg{\Delta}_i=(-1)\left|\begin{matrix}\*A_{11} & \cdots & \*A_{1,i-1} & \*A_{10} & \*A_{1,i+1} & \cdots & \*A_{1r}\cr
\vdots &  & \vdots & \vdots & \vdots & & \vdots \cr 
\*A_{r1} & \cdots & \*A_{r,i-1} & \*A_{r0} & \*A_{r,i+1} & \cdots & \*A_{rr}\end{matrix}\right|_\otimes,\qquad i=1,\ldots,r,
$$
are called operator determinants. For details see, e.g., \cite{AtkinsonBook}. Note that \eqref{eq:opdet} is a generalized Leibniz formula for the determinant,
where we sum over all permutations $\sigma$ in the permutation group $S_r$, with the Kronecker product instead of the 
usual product.
If $\vg{\Delta}_0$ is nonsingular, the problem \eqref{eq:add_mupa} is \emph{regular}, and the matrices 
$\vg{\Delta}_0^{-1}\vg{\Delta}_1,\ldots,\vg{\Delta}_0^{-1}\vg{\Delta}_r$ commute. A regular 
$r$-parameter eigenvalue problem \eqref{eq:add_mupa} has $\prod \limits_{i=1}^r n_i$ eigenvalues.

\section{Apply shift to singular problem}\label{appendix}
\noindent
The matrices $\vg\Delta_i$, $i=0,\ldots,4$, associated with the four-parameter eigenvalue problem \eqref{eq:4par} are such that $\vg\Delta_0$ is singular, but 
$\vg\Delta_0+s \vg\Delta_4$ is nonsingular for a generic shift $s\ne 0$. Thus, it is more efficient to solve a shifted system of generalized eigenvalue problems
\[
\vg\Delta_1 \*z = \I \widetilde k (\vg \Delta_0 + s \vg\Delta_4)\*z,\quad
\vg\Delta_2 \*z = \I \widetilde \kappa_{y,1} (\vg \Delta_0 + s \vg\Delta_4)\*z,\quad
\vg\Delta_3 \*z = \I \widetilde \kappa_{y,1} (\vg \Delta_0 + s \vg\Delta_4)\*z,\quad
\vg\Delta_4 \*z = \widetilde \xi_0 (\vg \Delta_0 + s \vg\Delta_4)\*z
\]
and then recover the eigenvalues of \eqref{eq:4par} as 
\[
\I k = \frac{\I \widetilde k}{1-s\widetilde \xi_0},\quad
\I \kappa_{y,1} = \frac{\I \widetilde \kappa_{y,1}}{1-s\widetilde \xi_0},\quad
\I \kappa_{y,2} = \frac{\I \widetilde \kappa_{y,2}}{1-s\widetilde \xi_0},\quad
\xi_0 = \frac{\widetilde \xi_0}{1-s\widetilde \xi_0}.\quad
\]
This approach circumvents the computationally complex staircase algorithm.
Note that the infinite eigenvalues of \eqref{eq:4par}  correspond to $1-s\widetilde\xi_0=0$.
The same idea can be applied to \eqref{eq:6par}, where $\vg\Delta_0+s \vg\Delta_6$ is 
nonsingular for a generic shift $s\ne 0$.

\section{Improve efficiency for isotropic plates}\label{app:isotropic}\noindent 
The quadratic eigenvalue problem \eqref{eq:EVPfreePlate} describing the free plate is typically solved at a given frequency by employing a linearization, resulting in a standard eigenvalue problem of twice the original size in $k$. If the plate consists of isotropic materials, there exists a transformation yielding a linear eigenvalue problem \textit{of the original size}. That is to say, it involves terms multiplied by $k^2$ but none in $k$ so that we can treat it as a linear eigenvalue problem in $\xi_0 = -k^2$. If the isotropic plate is coupled to unbounded fluid media, we can apply this approach analogously to our current formulation after minor modifications concerning the interface conditions. Consider the problem \eqref{eq:4par} and sort the degrees of freedom into horizontal displacements, vertical displacements, and pressures. As is known in the case of free plates \cite{Kausel1981,Gravenkamp2018}, the resulting finite-element matrices have a particular block structure:
\begin{multline}
    \left(-
    \left[\begin{array}{ccc}  
      \Etbar^{xx} &  \vt{0} & \vt{0}\\
      \vt{0} & \Etbar^{yy} & \Etbar^{yp}\\
      \vt{0} & \vt{0} & \vt{0}
    \end{array}\right]
    +\omega^2
    \left[\begin{array}{ccc}  
      \Mbar^{xx} &  \vt{0} & \vt{0}\\
      \vt{0} & \Mbar^{yy} & \vt{0}\\
      \vt{0} & \Mbar^{py} & \vt{0}
    \end{array}\right]
    + \I k
    \left[\begin{array}{ccc}  
      \vt{0} & \Eobar^{xy} & \vt{0}\\
      \Eobar^{yx} & \vt{0} & \vt{0}\\
      \vt{0} & \vt{0} & \vt{0}
    \end{array}\right]
    + (\I k)^2 
    \left[\begin{array}{ccc}  
     \Ezbar^{xx} &  \vt{0} & \vt{0}\\
     \vt{0} & \Ezbar^{yy} & \vt{0}\\
     \vt{0} & \vt{0} & \vt{0}
   \end{array}\right]\right. \\
   \left.
    + \I \kappa_{y,1}
    \left[\begin{array}{ccc}  
      \vt{0} & \vt{0} & \vt{0}\\
      \vt{0} & \vt{0} & \vt{0}\\
      \vt{0} & \vt{0} & \*R_1^{pp}
    \end{array}\right]
     + \I \kappa_{y,2} 
     \left[\begin{array}{ccc}  
      \vt{0} & \vt{0} & \vt{0}\\
      \vt{0} & \vt{0} & \vt{0}\\
      \vt{0} & \vt{0} & \*R_2^{pp}
    \end{array}\right]
    \right)
     \begin{bmatrix}
      \un^x \\ \un^y \\ \pn
    \end{bmatrix}
     =\vt{0}.   \label{eq:4par_blocks}
\end{multline}
We multiply the second and third equations in the above system by $\I k$ and modify the eigenvectors to obtain
\begin{multline}
  \left(-
  \left[\begin{array}{ccc}  
    \Etbar^{xx} &  -\Eobar^{xy} & \vt{0}\\
    \vt{0} & \Etbar^{yy} & \Etbar^{yp}\\
    \vt{0} & \vt{0} & \vt{0}
  \end{array}\right]
  +\omega^2
  \left[\begin{array}{ccc}  
    \Mbar^{xx} &  \vt{0} & \vt{0}\\
    \vt{0} & \Mbar^{yy} & \vt{0}\\
    \vt{0} & \Mbar^{py} & \vt{0}
  \end{array}\right]
  + (\I k)^2 
  \left[\begin{array}{ccc}  
   \Ezbar^{xx} &  \vt{0} & \vt{0}\\
   \Eobar^{yx} & \Ezbar^{yy} & \vt{0}\\
   \vt{0} & \vt{0} & \vt{0}
 \end{array}\right]\right. \\
 \left.
  + \I \kappa_{y,1}
  \left[\begin{array}{ccc}  
    \vt{0} & \vt{0} & \vt{0}\\
    \vt{0} & \vt{0} & \vt{0}\\
    \vt{0} & \vt{0} & \*R_1^{pp}
  \end{array}\right]
   + \I \kappa_{y,2} 
   \left[\begin{array}{ccc}  
    \vt{0} & \vt{0} & \vt{0}\\
    \vt{0} & \vt{0} & \vt{0}\\
    \vt{0} & \vt{0} & \*R_2^{pp}
  \end{array}\right]
  \right)
   \begin{bmatrix}
    \un^x \\ \I k\, \un^y \\ \I k\, \pn
  \end{bmatrix}
   =\vt{0}.   \label{eq:4par_blocksLinear}
\end{multline}
Hence, we have obtained a three-parameter eigenvalue problem with parameters $\I\kappa_{y,1},\,\I\kappa_{y,2},\,\xi_0$ that is solved analogously to \eqref{eq:4par} but without requiring \eqref{eq:4par_main_d}. Consequently, the operator determinants are of size $4n_1\times 4n_1$ (instead of $8n_1\times 8n_1$), leading to a less expensive assembly and solution of the generalized eigenproblem \eqref{eq:genEVP}.

\bibliographystyle{elsarticle-num}
\bibliography{embeddedWaveguide.bib, additional.bib}

\end{document}